\documentclass[12pt]{amsart}

\usepackage{fancyhdr}

\usepackage{amssymb,latexsym,amscd}

\usepackage{amsthm,amsmath}

\input xy

\xyoption{all}

\usepackage[cp850]{inputenc}

\newcommand{\gfr}{\mathfrak{g}}

\newcommand{\cO}{\mathcal{O}}

\newcommand{\cE}{\mathcal{E}}
\newcommand{\cG}{\mathcal{G}}

\newcommand{\cM}{\mathcal{M}}
\newcommand{\cA}{\mathcal{A}}

\newcommand{\cL}{\mathcal{L}}
\newcommand{\DD}{\mathcal{D}}

\newcommand{\Prym}{\mathrm{Prym}}

\newcommand{\Pic}{\mathrm{Pic}}
\newcommand{\Jac}{\mathrm{Jac}}
\newcommand{\Stab}{\mathrm{Stab}}
\newcommand{\Div}{\mathrm{Div}}
\newcommand{\Aut}{\mathrm{Aut}}
\newcommand{\Mor}{\mathrm{Mor}}

\newcommand{\lra}{\longrightarrow}

\newcommand{\hra}{\hookrightarrow}

\newcommand{\ra}{\rightarrow}

\newcommand{\la}{\lambda}

\newcommand{\vf}{\varphi}

\newcommand{\PP}{\mathbb{P}}
\newcommand{\Sbb}{\mathbb{S}}
\newcommand{\ZZ}{\mathbb{Z}}

\newcommand{\QQ}{\mathbb{Q}}
\newcommand{\NN}{\mathbb{N}}

\newcommand{\CC}{\mathbb{C}}

\newcommand{\End}{\mathrm{End}}

\newcommand{\Hom}{\mathrm{Hom}}
\newcommand{\Ind}{\mathrm{Ind}}
\newcommand{\Res}{\mathrm{Res}}

\newcommand{\GL}{\mathrm{GL}}
\newcommand{\SL}{\mathrm{SL}}

\newcommand{\Spin}{\mathrm{Spin}}
\newcommand{\Nm}{\mathrm{Nm}}

\newcommand{\im}{\mathrm{im}}
\newcommand{\coker}{\mathrm{coker}}
\def\map#1{\ \smash{\mathop{\longrightarrow}\limits^{#1}}\ }

\theoremstyle{plain}

\newtheorem{thm}{Theorem}[section]

\newtheorem{lem}[thm]{Lemma}

\newtheorem{prop}[thm]{Proposition}

\newtheorem{cor}[thm]{Corollary}

\newtheorem{rem}[thm]{Remark}
\newtheorem{defi}[thm]{Definition}

\begin{document}

\title[]{Polarizations of Prym varieties for Weyl groups via abelianization}

\author{Herbert Lange}
\author{Christian Pauly}
\address{Mathematisches Institut \\ Universit\"at Erlangen-N\"urnberg \\
Bismarckstrasse 1 1/2 \\ D-91054 Erlangen \\ Deutschland}

\email{lange@mi.uni-erlangen.de}
\address{D\'epartement de Math\'ematiques \\ Universit\'e de Montpellier II - Case Courrier 051 \\ Place Eug\`ene Bataillon \\ 34095 Montpellier Cedex 5 \\ France}
\email{pauly@math.univ-montp2.fr}

%\thanks{}

%\keywords{}

%\subjclass[2000]{Primary 14H60, 14D20, Secondary 14H40}

\begin{abstract}
Let $\pi: Z \ra X$ be a Galois covering of smooth projective curves with Galois group 
the Weyl group of a simple and simply connected Lie group $G$. For any dominant weight $\lambda$ consider the curve 
$Y = Z/\Stab(\lambda)$. The Kanev correspondence defines an abelian subvariety $P_\lambda$ of the Jacobian of $Y$. 
We compute the type of the polarization of the restriction of the canonical principal polarization of $\Jac(Y)$ to $P_\lambda$
in some cases. In particular, in the case of the group $E_8$ we obtain families of Prym-Tyurin varieties. The main idea is the use of 
an abelianization map of the Donagi-Prym variety to the moduli stack of principal $G$-bundles on the curve $X$.
 
\end{abstract}

\maketitle

%%%%%%%%%%%%%%%%%%%%%%%%%%%%%%%%%%%%%%%%%%%%%%%%%%%%%%%%%%

\section{Introduction}

\subsection{Verlinde spaces}

Let $X$ be a smooth complex projective curve of genus $g$ and let $G$ be a simple, simply-connected
complex Lie group. We denote by $\cM_X(G)$ the moduli stack of principal $G$-bundles and by 
$\cL$ the ample generator of its Picard group. The celebrated Verlinde formula (\cite{Fa1},
\cite{So1}, \cite{S}) computes the dimension $N_{g,l}(G)$ of the space of global sections
$H^0(\cM_X(G), \cL^{\otimes l}) $  for any level $l$. The Verlinde numbers at level
$l=1$ for the groups of type ADE are given in the following table.

\bigskip

\begin{center}
\begin{tabular}{|c||c|c|c|c|c|}
\hline 
 $G$ & $\SL(m)$ & $\Spin(2m)$ & $E_6$ & $E_7$ & $E_8$  \\
\hline
 $N_{g,1}(G)$ & $m^g$ & $4^g$ & $3^g$ & $2^g$ & 1 \\
\hline 
\end{tabular}
\end{center}

\bigskip

The number $m^g$ for $\SL(m)$ coincides with the number of level-$m$ theta functions on the Jacobian of 
$X$ (see \cite{BNR}). For the even Spin group the Verlinde number equals the number of theta characteristics
of $X$ (see \cite{O}). The striking simplicity of the Verlinde numbers for $E_6$, $E_7$ and
$E_8$ was the main motivation for us to try to relate these Verlinde spaces to spaces of
theta functions on polarized abelian varieties (the Prym varieties) and compute the
induced polarizations --- see Main Theorem.

\subsection{Abelianization of principal $G$-bundles}

The abelianization program of principal $G$-bundles, or more precisely $G$-Higgs bundles, takes its
origin in Hitchin's papers \cite{Hi1} and \cite{Hi2}. For the case $G = \SL(m)$ it is shown in \cite{BNR}
that for a sufficiently ramified spectral cover $\psi: Y \ra X$ the direct image map
$$ \Prym(Y/X) \lra \cM_X(\SL(m)) $$
induces by pull-back an isomorphism between the $\SL(m)$-Verlinde space at level $1$ and the space
of abelian theta functions $H^0(\Prym(Y/X), L_Y)$. For general structure groups $G$ the abelianization
theory has been worked out by  Faltings \cite{Fa2}, by  Donagi \cite{D}, \cite{D2} and Donagi-Gaitsgory
\cite{DG}.

\subsection{Correspondences on spectral and cameral covers}

For general structure groups $G$ Prym varieties can be constructed via correspondences on covers of the 
curve $X$:

\bigskip

In \cite{K} Kanev constructs from the data of a rational map $f : \CC \ra \gfr = \mathrm{Lie}(G)$ and an
irreducible representation $\rho_\lambda  : G \ra \GL(V)$ a spectral cover $\psi: Y \ra \PP^1$ equipped with a
correspondence. He shows that if $G$ is of type ADE, the weight $\lambda$ minuscule and $f$ sufficiently general, then
the Prym variety  $P_\lambda \subset \Jac(Y)$ associated with Kanev's correspondence is a Prym-Tyurin
variety (i.e. the polarization on $P_\lambda$ induced from the principal polarization of $\Jac(Y)$ is a
multiple of a principal polarization). Note that Kanev's construction is carried out in the case $X = \PP^1$.

\bigskip

A different but related construction of Prym varieties is given by Donagi in \cite{D}, \cite{D2}:
let $T \subset G$ be a maximal torus, $W$ the Weyl
group of $G$ and $\Sbb_{\omega} = \Hom(T,\CC^*)$ the weight lattice. For any cameral cover, i.e. a Galois
cover with Galois group $W$
\[
\pi:Z \lra X
\]
satisfying some conditions on the ramification, Donagi introduces the Prym variety 
$\Prym(\pi,\Sbb_{\omega}) := \Hom_W(\Sbb_{\omega},\Jac(Z))$ parametrizing $W$-equivariant
homomorphisms from $\Sbb_{\omega}$ to $\Jac(Z)$.

\bigskip

In this paper we generalize Kanev's construction to an arbitrary base curve $X$. Given a Galois
cover $\pi : Z \ra X$ and a  dominant weight 
$\lambda \in \Sbb_{\omega}$ we consider the cover of curves
$$ \psi: Y \lra X, \qquad \text{with} \qquad 
Y = Z/\Stab(\lambda). $$
Kanev's construction generalizes (see section 3) to give a correspondence $\overline{K}_{\lambda}$ 
on the curve $Y$ defining an
abelian subvariety
\[
P_{\lambda} \subset \Jac(Y),
\]
which is isogenous to the Donagi-Prym variety $\Prym(\pi,\Sbb_{\omega})$ (Proposition 6.13).

\subsection{The polarization on the Prym variety $P_\lambda$}

Let $L_Y$ denote a line bundle defining the canonical principal polarization on the 
Jacobian $\Jac(Y)$. The aim of this paper is to 
compute the induced polarization $L_Y|P_{\lambda}$ under certain assumptions. In fact, if $q_{\lambda}$ denotes the exponent 
of the correspondence $\overline{K}_{\lambda}$ and $d_{\lambda}$ the Dynkin index of $\lambda$, our main result is 
the following theorem (we prove a slighty more precise version, see Theorem \ref{maintheo}) . We use the notation of \cite{Bou} for the weights.

\bigskip
\noindent
{\bf Main Theorem.} {\em
We suppose that the $W$-Galois cover $\pi: Z \ra X$ is \'etale. Then in the cases given in the table below the 
induced polarization $L_Y|P_\lambda$ is divisible by $q_\lambda$,i.e. 
$L_Y|P_\lambda = M^{\otimes q_\lambda}$ and the polarization $M$ on $P_\lambda$ is of type 
$K(M) = (\ZZ/m\ZZ)^{2g}$: }

\bigskip

\begin{center}
\begin{tabular}{|c|c|c||c|}
\hline 
 Weyl group of type & weight $\la$ & $q_\la = d_\la$ & $K(M)$ \\ \hline \hline
  $A_n$ & $\varpi_i ; \ (i,n+1) = 1$ & ${n-1 \choose i-1}$   & $(\ZZ/(n+1)\ZZ)^{2g}$ \\ \hline 
$D_n$; \ $n$ odd & $\varpi_{n-1}, \varpi_n$  & $2^{n-3}$ & $(\ZZ/4\ZZ)^{2g}$ \\ \hline
$E_6$ & $\varpi_1,\varpi_6$ &  $6$ & $(\ZZ/3\ZZ)^{2g}$ \\ \hline
$E_7$ & $\varpi_7$ & $12$ &$(\ZZ/2\ZZ)^{2g}$\\ \hline
$E_8$ & $\varpi_8$ & $60$ & $0$ \\ \hline 
\end{tabular}
\end{center}

\bigskip

We observe that in all the cases of this table we have an equality (see table in section 1.1)
$$ \dim H^0(P_\la,M) = N_{g,1}(G).$$ 
A priori it is not clear (e.g. in the case $G= E_8$) whether the linear map between the
Verlinde space and
$H^0(P_\la,M)$ induced by the abelianization map $\Delta_\theta$ ---
see below --- is nonzero.

\bigskip

It is well known that the  Weyl group of type $E_k$  is closely related to the del Pezzo surface of degree $9-k$ 
for $5 \leq k \leq 8$. In fact, a slightly modified lattice of the weight lattice is isomorphic to the 
Picard lattice of the corresponding del Pezzo surface (see \cite[section 8.7]{K}). 
Moreover, for $4 \leq k \leq 7$ the Kanev correspondence is given essentially by the incidence correspondence of
 lines of the corresponding del Pezzo surface. For $k=8$ there are multiplicities due to the fact that 
 the weight $\varpi_8$ is only quasi-minuscule (see \cite{K}).    
Notice that in these cases the polarization $M$ on $P_\lambda$ is of type $(\ZZ/d\ZZ)^{2g}$ where $d$ is the degree of 
the correponding del Pezzo surface. 

\bigskip

In particular in the case of $W = W(E_8)$ 
we obtain a family of Prym-Tyurin varieties, i.e. the pairs $(P_{\lambda},M)$ are principally polarized abelian varieties.
It is easy to see that any curve $X$ of genus $g \geq 4$ admits an \'etale Galois covering with Galois group $W(E_8)$ 
and we get a family of Prym-Tyurin varieties of dimension $8(g-1)$ of exponent 60. We plan to study this family in a 
subsequent paper.\\

Probably there is an analogous result in the case that the Galois covering $\pi: Z \ra X$ admits simple ramification. Certainly
the paper \cite{DG} will be essential for this. We plan to come back to it subsequently. \\

Note that our results are disjoint from the results in \cite{K}. First of all, Kanev considers only Galois coverings over 
$\PP^1$ which are necessarily ramified. Moreover, his correspondence satisfies a quadratic equation and he uses his criterion 
(\cite{K2}) to show that the associated abelian subvarieties are principally polarized. In our case the corresponding 
correspondence satisfies a cubic equation (see Theorem \ref{thm3.9}) and a quadratic equation only on the Prym variety $\Prym(Y/X)$
of the covering $\psi: Y \ra X$. However, the abelian variety $\Prym(Y/X)$ is not principally polarized and thus we cannot apply 
Kanev's criterion \cite{K2} in order to compute the polarization of $P_\lambda$.

Instead we proceed as follows: We work out a general result on restrictions of polarizations to abelian subvarieties 
(Proposition \ref{prop2.10}) roughly saying that, if the restricted polarization equals 
the $q$-fold of a polarization, where 
$q$ is the exponent of the abelian subvariety, then its type can be computed.

The main idea of the proof is then the use of an abelianization map 
\[
\Delta_{\theta} : \Prym(\pi,\Sbb_{\omega})_n \ra \cM_X(G),
\] 
(see Section 7). Here $\Prym(\pi,\Sbb_{\omega})_n$ denotes a certain component of the Donagi-Prym variety 
$\Prym(\pi,\Sbb_{\omega})$. 
The fact, that the restricted polarization is the $q_\lambda$-fold of a polarization 
is a consequence of the existence of a commutative diagram
$$
\begin{CD}
\Prym(\pi, \Sbb_\omega)_n @>\widetilde{ev}_\lambda>> T_\alpha(P_\lambda )  \\
@VV \Delta_\theta V @VV \psi_* V \\
\cM_X(G) @>\tilde{\rho}_\lambda>> \cM_X(\mathrm{SL}(m))
\end{CD}
$$
(here $\widetilde{ev}_\lambda$ denotes evaluation at $\lambda$, the map $\tilde{\rho}_\lambda$ is induced by the representation with  dominant weight $\la$, 
and $T_\alpha$ is a certain translation)
together with a theorem of Laszlo-Sorger (\cite{LS}, \cite{S}) saying that the pull-back of the determinant bundle 
on $\cM_X(\mathrm{SL}(m))$ by $\tilde{\rho}_\lambda$ is the $d_\lambda$-fold of the ample generator of $\Pic(\cM_X(G))$ and 
the assumption $q_\lambda = d_\lambda$. \\

The contents of the paper is as follows:
In Section 2 we prove the result just mentioned concerning restriction of  a polarization to an abelian subvariety. 
In Section 3 we prove that the Schur- and the Kanev correspondence belong to the same abelian subvariety.
This was shown in \cite{LRo} for the special case $X = \PP^1$ and for a left action of the group. In Section 4 we compute
the invariants of our main examples mentioned above. In Section 5 we introduce the Donagi-Prym variety and derive the properties 
we need. Section 6 contains some results on Mumford groups. In Section 7 we introduce the abelianization map and prove
the commutativity of the above diagram. Finally, Section 8 contains the proof of our main theorem. \\

Finally one word to the group actions: The group $W$ acts on the curve $Z$ as well as on the weight lattice $\Sbb_\omega$
and these actions have to be consistently either both left actions or both right actions. Since traditionally principal 
bundles are defined via right actions, we are forced to use right actions in both cases.\\

We would like to thank Laurent Manivel for a helpful discussion.

\bigskip

\section{Restriction of polarizations to abelian subvarieties}

\noindent
Let $S$ be an abelian variety with polarization $L$ and associated isogeny $\vf_L: S \lra \hat{S}$. Denote as usual
$K(L) = \ker \vf_L$. We consider an endomorphism $u \in \End(S)$ satisfying the following
properties
\begin{itemize}
\item $u$ is symmetric with respect to the Rosati involution, i.e., 
\begin{equation}\label{Rosati}
\varphi_L \, u = \hat{u} \, \varphi_L,
\end{equation}
where $\hat{u} \in \End(\hat{S})$ is the dual endomorphism.
\item there exists a positive integer $q$ such that 
\begin{equation} \label{quadeq}
u^2 = q u.
\end{equation}
\end{itemize}
We introduce the abelian subvarieties of $S$ 
\begin{equation} \label{defprym}
A = \im (u) \subset \ker(q-u) \qquad \text{and} \qquad  P = \im (q-u) \subset \ker(u)
\end{equation}
and we wish to study the induced polarizations $L_A = L|A$ on $A$ and $L_P = L|P$ on $P$, i.e. determine the subgroups 
$K(L_P)$ and $K(L_A)$.

In the sequel we use the following notation
$$ \iota_A: A \hra S, \qquad \iota_P: P \hra S, \qquad \pi_A: S \lra S/A, \qquad  \pi_P: S \lra S/P.$$
where the $\iota$'s are inclusions and $\pi$'s are projections. 

\begin{rem} \label{remarkprym}
{\em (i): For any abelian subvariety $A$ of $S$ there is an endomorphism $u$ satisfying \eqref{Rosati} and \eqref{quadeq}. In fact, 
this is valid for the composition 
$$
u = \iota_A \, \psi_{L_A} \, \hat{\iota}_A \, \vf_L \in \End (S)
$$ 
(see \cite[Lemma 5.3.1]{BL}). Here $\psi_{L_A} = q \, \vf_L^{-1}: \hat{A} \ra A$ is an isogeny.

(ii):  Let $\psi: Y \lra X$ be a degree $d$ cover of smooth projective curves. If we identify the Jacobians $\Jac(X)$ and $\Jac(Y)$
with their duals via the canonical polarizations, then the dual
$\hat{\psi}^*$ of the map $\psi^*: \Jac(X) \lra \Jac(Y)$ coincides with the norm map 
$\Nm: \Jac(Y) \lra \Jac(X)$. Moreover the endomorphism
$$
t = \psi^* \,  \Nm \in \End(\Jac(Y))
$$
satisfies \eqref{Rosati} and \eqref{quadeq} with $q = d$. 
The abelian subvariety $P = \im(d-t)$ defined in \eqref{defprym} coincides with the usual Prym variety (not necessarily principally polarized) of the morphism $\psi$, which we denote by $\Prym(Y/X)$, i.e.,
$$ P = \Prym(Y/X) = (\ker \Nm)_0.$$
The endomorphism $t$ is induced by the trace correspondence
$\overline{T}$ on $Y$ --- see section 3.
We denote the kernel $\ker \psi^* : \Jac(X) \lra \Jac(Y)$ by $K$. We also have the
equality 
$$ \im (\psi^*) = \Jac(X)/K = \im (t). $$
Since the norm
map $\Nm$ is the dual of $\psi$ we easily obtain that the group of connected components
of the fibre $\Nm^{-1}(0)$ is isomorphic to $K$.

(iii): With $u$ also $u' = q -u$ satisfies \eqref{Rosati} and \eqref{quadeq}. Replacing $u$ by $u'$ interchanges the 
abelian subvarieties
$A$ and $P$. Hence for the following Lemmas it suffices to prove one of the two statements concerning $A$ and $P$ and we will 
do so without further mentioning.
}
\end{rem}

\bigskip

Recall that the addition map $\mu$ of $S$ satisfies 
$\mu = \iota_A + \iota_P$ 
and that the following sequence is exact 

\begin{equation} \label{eq3}0 \lra  A \cap P \lra A \times P \map{\mu} S \lra 0.
\end{equation}
For any abelian variety $B$ and any positive integer $n$ let $B_n$ 
denote the subgroup of $n$-torsion points of $B$. Then we have, denoting by $|M|$ the cardinality of a set $M$,

\begin{lem} \label{lem2.2}
The finite subgroup $A \cap P$ is contained in $A_q$ and $P_q$.
\end{lem}

\begin{proof}
Let $x \in A \cap P$. Then $u(x)=qx$, since $x \in A$ and $u(x) = 0$, since $x \in P$. This implies the assertion.
\end{proof}

\begin{lem} \hspace{0.5cm} \label{lem2.3}
$|K(L_A)|\cdot |K(L_P)| = |A \cap P|^2 \cdot |K(L)|.$
\end{lem}

\begin{proof}
According to \cite[Corollary 5.3.6]{BL} the polarization $\mu^*(L)$ splits, i.e. 
$\vf_{\mu^*(L)} = \vf_{L_A} \times \vf_{L_P}$, which gives
$$
|K(\mu^*L)| = |K(L_A)| \cdot |K(L_P)|.
$$
On the other hand, (\ref{eq3}) implies $\deg(\mu) = |A \cap P|$ and thus
$$
|K(\mu^*L)| = |\ker(\mu\vf_L\hat{\mu})| = \deg(\mu)^2 \cdot |\ker(\vf_L)| = |A \cap P|^2 \cdot |K(L)|.
$$
Combining both equations gives the assertion.
\end{proof}

\begin{lem} \label{lem2.4}
We have the following equalities of abelian subvarieties of $\hat{S}$,
$$ \varphi_L(P) = \widehat{S/A} \subset \hat{S} \qquad \text{and} \qquad \varphi_L(A) = \widehat{S/P} \subset \hat{S}.$$
\end{lem}

\begin{proof}
The endomorphism $u$ factorizes as $i_A \, v$ with $v: S \ra A$. 
Taking the dual gives the factorization
$$ \hat{u} : \hat{S} \map{\hat{\iota}_A} \hat{A} \map{\hat{v}} \hat{S},$$
from which we deduce that 
$$ \ker \hat{\iota}_A = \widehat{S/A} \subset \ker \hat{u}.$$
On the other hand we have $P \subset \ker u$ and so \eqref{Rosati} implies
$\varphi_L(P) \subset \ker \hat{u}$. This gives the
first equality, since both abelian subvarieties are the connected component of the origin of $\ker \, \hat{u}$, as they are 
of the same dimension.  
\end{proof}

\bigskip

We will give two descriptions of the subgroups 
$K(L_P) \; (= \ker (\varphi_{L_P} :  P \map{\iota_P} S \map{\varphi_L} \hat{S} \map{\hat{\iota}_P} \hat{P} ))$ and $K(L_A)$.
We observe that $\ker \hat{\iota}_P = \widehat{S/P}$ and, 
according to Lemma \ref{lem2.4}, $ \widehat{S/P} = \varphi_L(A)$. Hence we obtain
$$ K(L_P) = \bigcup_{x\in K(L)} (A+x) \cap P,$$
where $A+x$ denotes the image of $A$ under translation by $x$. Similarly we have
$$ K(L_A) = \bigcup_{x\in K(L)} (P+x) \cap A.$$
In particular,
$$
A \cap P \subset K(L_P) \quad \mbox{and} \quad A \cap P \subset K(L_A).
$$\\
For the second description consider the isogenies
$$ \alpha = \pi_A \, \iota_P : P \map{\iota_P} S \map{\pi_A} S/A \qquad \text{with} \qquad
\ker \alpha = A \cap P$$
and
$$\beta = \pi_P \, \iota_A : A \map{\iota_A} S \map{\pi_P} S/P \qquad \text{with} \qquad
\ker \beta = A \cap P.$$
Moreover, define the isogenies 
$$
\vf_P : P \lra \widehat{S/A} \quad \mbox{and} \quad \vf_A: A \lra \widehat{S/P}
$$ by the factorizations $\vf_L|P: P \map{\vf_P} \widehat{S/A} \hra \hat{S}$ and 
$\vf_L|A: A \map{\vf_A} \widehat{S/P} \hra \hat{S}$ of Lemma \ref{lem2.4}. With this notation we have

\begin{lem} \label{lem2.5}
The following sequences are exact
\begin{equation} \label{eq4}
0 \lra A \cap P \hra K(L_P) \map{\alpha} \ker (\hat{\vf}_P) \lra 0 \quad \text{with} \quad
K(L) \cap P \cong \ker(\hat{\vf}_P) \subset S/A,
\end{equation}
\begin{equation} \label{eq5}
0 \lra A \cap P \hra K(L_A) \map{\beta} \ker (\hat{\vf}_A) \lra 0 \quad \text{with} \quad
K(L) \cap A \cong \ker (\hat{\vf}_A) \subset S/P.
 \end{equation}
\end{lem}

\begin{proof}
The dual isogeny of $\alpha$ factorizes as
$$ \hat{\alpha}: \widehat{S/A} \map{\hat{\pi}_A} \hat{S} \map{\hat{\iota}_P} 
\hat{P}.
$$
By the previous lemma $\varphi_L(P) = \widehat{S/A}$, which implies that
$$ \vf_{L_P} = \hat{\iota}_P \, \varphi_L \, \iota_P = \hat{\alpha} \, \varphi_L \, \iota_P.$$
Hence $\varphi_{L_P}$ factorizes as  $\varphi_{L_P} : P \map{\vf_P} \widehat{S/A} \map{\hat{\alpha}} \hat{P}.$
Taking the dual and using $\hat{\varphi}_{L_P} = \vf_{L_P}$ we obtain the following factorization
$$ \varphi_{L_P}: P \map{\alpha} S/A \map{\hat{\vf}_P} \hat{P}$$
which gives the exact sequence \eqref{eq4}. The assertion on the kernel follows from 
$\ker (\hat{\vf}_P) \cong \ker (\vf_P) = K(L) \cap P$.
\end{proof}

In particular we obtain
$$ |K(L_P)| = |A \cap P|\cdot |K(L) \cap P| \quad \mbox{and} \quad |K(L_A)| = |A \cap P|\cdot |K(L) \cap A|.$$
Inserting this into Lemma \ref{lem2.3} we conclude

\begin{cor} \label{cor2.5} \hspace{0.5cm} $
|K(L) \cap A| \cdot |K(L) \cap P| = |K(L)|.
$
\end{cor}

\bigskip

Consider the finite groups
$$
G_P = \ker u /P \qquad \mbox{and} \qquad G_A = \ker (q-u) /A.
$$
They are the groups of connected
components of $\ker u$ and $\ker (q-u)$ respectively. 

\begin{lem} \label{lem2.7}
There are canonical exact sequences
$$ 0 \lra A \cap P \lra A_q \lra G_P \lra 0 \quad \mbox{and} \quad  0 \lra A \cap P \lra P_q \lra G_A \lra 0 .$$
\end{lem}

\begin{proof}

Since $P \subset \ker u$ the endomorphism $u$ descends to an isogeny
$$ \bar{u}: S/P \lra A, \qquad \text{with} \qquad \ker \bar{u} = G_P. $$
Moreover we know that $u|A = q$, which implies that the
composite isogeny
$$ A \map{\pi_P \, \iota_A} S/P \map{\bar{u}} A$$
equals $q$. Taking kernels leads to the first exact sequence of the Lemma.
\end{proof}

We denote by 
$$
u_q : S_q \lra A_q \quad \mbox{and} \quad (q-u)_q : S_q \lra P_q
$$ the restrictions of $u$  and $q-u$ to the $q$-torsion
points $S_q$. Note that $u_q$ and $(q-u)_q$ are nilpotent:  $u_q^2 = (q-u)_q^2 = 0$. Note moreover that 
$G_P \subset \left( S/P \right)_q$ and $G_A \subset \left( S/A \right)_q$.

\begin{lem}
We have the following equalities of finite groups
$$ \im \, u_q = A \cap P, \quad \coker \, u_q = G_P \quad \mbox{and} \quad \im \, (q-u)_q = A \cap P, \quad  \coker \, (q-u)_q = G_A.$$
\end{lem}

\begin{proof}
First we claim that $\im \, u_q \subset A \cap P$. For the proof let $a \in \im \, u_q$. Write $a= u(s)$ with
$s \in S_q$. Then $qs=0$ and therefore $a= (u-q)(s) \in P$. 

In order to show equality of these subgroups of $A_q$ we compute their indices. According to Lemma \ref{lem2.7} 
the index of $A\cap P$ in $A_q$ is $|G_P|$. The index of $\im \, u_q$ in $A_q$ equals 
$$ \frac{|A_q|}{|\im \, u_q|} = \frac{|A_q|}{|S_q|} \cdot |\ker u_q| = \frac{|\ker u_q|}{|P_q|}.$$
Moreover since $G_P \subset \left( S/P \right)_q$ the exact sequence 
$$ 
0 \lra P \lra \ker u \lra G_P \lra 0$$
remains exact after taking $q$-torsion points. This implies the first equality. The second equality is obvious.
\end{proof}

Consider again the isogenies $\hat{\vf}_P: S/A \lra \hat{P}$ and 
$\hat{\vf}_A: S/P \lra \hat{A}$.

\begin{lem} \label{lem2.9}
We have the following equalities of subgroups of $S/A$ and $S/P$ respectively:
$$ \ker \hat{\vf}_P = \pi_A (K(L)) \quad \mbox{and} \quad \ker \hat{\vf}_A = \pi_P(K(L)). $$
\end{lem}

\begin{proof}
The inclusion $\pi_A(K(L)) \subset \ker \hat{\vf}_P$ follows from
the commutative diagram
$$
\begin{CD}
S @>\varphi_L>> \hat{S} \\
@VV\pi_AV    @VV\hat{i}_PV \\
S/A @>\hat{\vf}_P>> \hat{P}
\end{CD}
$$
So it will be enough to show that
these two groups have the same order. We have, using \eqref{eq4},
$$|\ker \hat{\vf}_P| = |K(L) \cap P| \qquad  \text{and} \qquad  
|\pi_A(K(L))| =  \frac{|K(L)|}{|K(L) \cap A|}.$$
Equality follows from Corollary \ref{cor2.5}.
\end{proof}

The following proposition will be applied in the proof of the main theorem (section 8).

\begin{prop}  \label{prop2.10}
{\em (a)}: If $L_A = M^q$ for a polarization $M$ on $A$, we have the following equality of subgroups of $A$,
$$ K(M) = u(K(L)),$$
{\em (b)}: If $L_P = N^q$ for a polarization $N$ on $P$, we have the following equality of subgroups of $P$,
$$ 
K(N) = (q-u)(K(L)).
$$
\end{prop}

Note that $L_A = M^q$ for a polarization $M$ on $A$ if and only if the polarization $L_A$ is divisible by $q$, 
meaning that $A_q \subset K(L_A)$. In this case 
$$K(M) = K(L_A)/A_q.$$ 

\begin{proof}
We take the quotient by $A \cap P$ of the inclusion $A_q \subset K(L_A)$. Since by Lemma \ref{lem2.7}, $A_q/A \cap P \cong G_P$  
and by Lemma \ref{lem2.5}, $K(L_A)/A \cap P \cong \ker \hat{\vf}_A$, this gives an inclusion
$$
G_P \subset \ker \hat{\vf}_A.
$$  
Recall that $G_P = \ker u/P$ and
by Lemma \ref{lem2.9}, $\ker \hat{\vf}_A = \pi_P (K(L))$. Therefore
$$
\begin{array}{ll}
 K(M) =  K(L_A)/A_q = \ker \hat{\vf}_A / G_P &= ((K(L) + P)/P)/(\ker u/P) \\
 &=  (K(L) + P)/\ker u =  u(K(L)).
 \end{array}
$$
\end{proof}

\section{The Schur and Kanev correspondences}

\subsection{Representations of the Weyl group $W$}

Let $W$ be a Weyl group and let $\widehat{W}$ denote the set of its irreducible characters.
It is known (see e.g \cite{Spr} Corollary 1.15) that all irreducible representations of $W$ 
are defined over $\QQ$. Therefore any irreducible representation of $W$ is also absolutely irreducible.  
Given $\omega \in \widehat{W}$  we choose an irreducible $\ZZ[W]$-module
$\mathbb{S}_{\omega}$ such that 
$$
V_\omega := \mathbb{S}_{\omega} \otimes_\ZZ \QQ
$$ 
is the irreducible 
representation  of $W$ corresponding to $\omega$.  As outlined in the introduction, 
we consider every representation of $W$ as a right-representation. In particular $\mathbb{S}_
{\omega}$ is a right-$\ZZ[W]$-module with action
$(\la,g) \mapsto \la g$ for $\la \in \mathbb{S}_{\omega}$ and $g \in W$.

\bigskip

Fix a {\it weight} $\la$
of the lattice $\mathbb{S}_{\omega}$, which is by definition (see \cite{K} page 158) a vector 
$$\la \in V_{\omega} \ \ \text{such that} \ \ 
\la g - \la \in \mathbb{S}_{\omega} \ \ \forall g \in W. $$

Let $\pi: Z \ra X$ denote a Galois covering of smooth projective curves with Galois group $W$. We consider the action of $W$ on $Z$ as a right-action $(z,g) \mapsto z^g$ 
for $z \in Z$ and $g \in W$. 
Let $H := \Stab(\la) \subset W$ denote the stabilizer subgroup of the weight $\lambda \in V_\omega$.
Then $\pi$ factorizes as follows
$$
\xymatrix{ Z \ar[dd]_{\pi}
\ar[dr]^{\varphi} \\ & Y \ar[dl]^{\psi} \\ X }
$$
where $Y$ denotes the quotient $Z/H$. Schur's orthogonality relations induce a correspondence on $Z$ which we denote by 
$S_\la$. On the other hand,
Kanev \cite{K} defined in the case $X = \PP^1$ a correspondence on the curve $Y$, which we denote by 
$\overline{K}_\la$. It is the aim of this 
section to generalize Kanev's construction to an arbitrary base curve $X$ and to work out the relation between the
two  correspondences $S_\la$ and $\overline{K}_\la$.

\bigskip

\subsection{Schur correspondence}

Since $V_\omega$ is an absolutely irreducible $W$-representation, there is a unique negative definite $W$-invariant symmetric form 
$(\; , \,)$ on $V_\omega$ with
\begin{enumerate}
\item $(\la, \mu) \in \ZZ$ for all $\mu \in \mathbb{S}_{\omega}$, 
\item any $W$-invariant form on $V_\omega$ satisfying (1) is an integer multiple of $(\; , \,)$.
\end{enumerate}
The {\it Schur correspondence} associated to the pair $(\mathbb{S}_{\omega}, \la)$ is by definition the rational correspondence on
$Z$ over $X$ defined by
$$
S_\la = \sum_{g \in W} (\la g, \la) \Gamma_g.
$$
Here $\Gamma_g \in Z \times Z$ denotes the graph of the automorphism $g$ of $Z$. Note that 
$( \la g, \la)$ need not be an integer.   
Considered as a map of $Z$ into the group 
$\Div_{\QQ} (Z)$ of rational divisors on $Z$, $S_\la$ is given by
$$
S_\la(z) = \sum_{g \in W} (\la g, \la)z^g.
$$

The correspondence $S_\la$ descends to a correspondence $\overline{S}_\la$ on $Y$ in the usual way, namely
$$
\overline{S}_\la = (\vf \times \vf)_*S_\la \subset Y \times Y.
$$
In order to express this correspondence as a map $Y \lra \Div_{\QQ}(Y)$, denote $d = [W:H]$ and 
let $\{g_1 = 1,g_2, \ldots, g_d\}$ denote a set of representatives for the right cosets of $H$ in $W$:
$$
W = \cup_{i=1}^d Hg_i.
$$

\begin{prop} \label{prop3.1}
For any $y = \vf (z) \in Y$ we have
$$
\overline{S}_\la(y) = |H|^2 \sum_{i=1}^d( \la g_i, \la) \vf(z^{g_i^{-1}}).
$$
\end{prop}

Note that the left-hand side of this equality does not depend neither on the choice of the 
point $z$ in the fibre over $y$ nor on the set of representatives $\{ g_i \}$.

\begin{proof} Since $H$ is the stabilizer of $\la$, we have by definition of $\overline{S}_\la$,
\begin{align*}
\overline{S}_\la(y)  & = \sum_{i=1}^d\sum_{h\in H}(\la hg_i, \la)\sum_{h' \in H}\varphi(z^{h'hg_i})
 = |H|\sum_{i=1}^d\sum_{k\in H}(\la g_i, \la )\varphi(z^{kg_i})  
  = |H| \sum_{g \in W} (\la g, \la) \varphi(z^g).
\end{align*}
Now, if $\{g_i\}$ is a set of representatives for the right cosets of $H$, then $\{g_i^{-1}\}$
is a set of representatives for the left cosets of $H$. Hence for any pair $(k,i) \in H \times \{1, \dots, d\}$
there is a unique pair $(h,j) \in H \times \{1, \dots, d\}$ such that $kg_i= g_j^{-1}h$.
Moreover, if $kg_i$ runs exactly once through $W$, so do the elements $g_j^{-1}h$. This implies, since
$(\;, \, )$ is $W$-invariant and since $H$ stabilizes $\la$,
$$
\overline{S}_\la(y) = |H|\sum_{i=1}^d\sum_{h \in H}( \la g_i, \la)\varphi(z^{g_i^{-1}h}) 
 = |H|\sum_{i=1}^d\sum_{h\in H}(\la g_i, \la) \varphi(z^{g_i^{-1}}).
$$
This implies the assertion.\end{proof}

\bigskip

\subsection{Kanev correspondence}

In order to define the Kanev correspondence on $Y$, let $U \subset X$ denote the complement of the branch locus of $\pi$ and
fix a point $\xi_0 \in U$. Since $H$ is the stabilizer of $\la$, the group $W$ acts on the set 
$\{ \la = \la g_1, \la g_2, \ldots, \la g_d\}$.

The group $W$ also acts on the fibre $\psi^{-1}(\xi_0)$. Note that the group $W$ does not act on the curve $Y$, since the 
subgroup $H$ is not necessarily a normal subgroup of $W$. However we can define the action of
$W$ on the set $\psi^{-1}(\xi_0)$ via the monodromy at the point $\xi_0$ as follows: consider the
fundamental group $\pi_1(U,\xi_0)$ and let $U_Z = \pi^{-1}(U) \subset Z$ and 
$U_Y = \psi^{-1}(U) \subset Y$. Choose a point $z \in Z$ such that $\pi(z) = \xi_0$. Then  
$\pi_1(U_Z,z)$
is a normal subgroup of $\pi_1(U,\xi_0)$ and we have the following
exact sequence 
$$ 0 \lra \pi_1(U_Z,z) \lra \pi_1(U,\xi_0) \lra W \lra 0.$$
The monodromy of the cover $\psi: U_Y \lra U$ at the point $\xi_0$ gives a homomorphism
$$\rho : \pi_1(U,\xi_0) \lra \mathrm{Aut}(\psi^{-1}(\xi_0)).$$
The kernel $\ker \rho$ equals $\cap_y \pi_1( U_Y ,y)$ where $y$ varies over the set
$\psi^{-1}(\xi_0)$. We note that $\pi_1(U_Z,z) \subset \pi_1( U_Y,y)$ for 
any point $y \in \psi^{-1}(\xi_0)$,  hence $\pi_1(U_Z,z) \subset \ker \rho$.
Therefore the monodromy map $\rho$ factorizes over $W$ 
$$\rho : \pi_1(U,\xi_0) \lra W \lra \mathrm{Aut}(\psi^{-1}(\xi_0)).$$
%Let $i$ denote the map $W \lra \mathrm{Aut}(\psi^{-1}(\xi_0))$. Then $\ker i$ 
%equals $\cap_w w^{-1} H w$ --- here $w$ varies over $W$ ---  which equals the kernel of the
%$W \lra GL(V)$. Do we assume that this representation is faithful?
Notice that according to our definitions the monodromy action of $\pi_1(U,\xi_0)$ on the fibre $\psi^{-1}(\xi_0)$
is a right action.\\

Choosing an element in the fibre $\psi^{-1}(\xi_0)$ induces a bijection
$$
\{ \la g_1, \ldots, \la g_d\} \stackrel{\sim}{\ra} \psi^{-1}(\xi_0),
$$
which is $W$-equivariant according to the definitions. In the sequel we identify the above sets,
 i.e. we label the elements of
$\psi^{-1}(\xi_0)$ by $\la = \la g_1, \ldots, \la g_d$.

For every point $\xi \in U$ choose a path $\gamma_{\xi}$ in $U$ connecting $\xi$ and $\xi_0$. The path
defines a bijection
$$
\mu: \psi^{-1}(\xi) \to \psi^{-1}(\xi_0)=\{\la  = \la g_1, \dots, \la g_d\}
$$
in the following way: For any $y \in \psi^{-1}(\xi)$ denote by $\tilde{\gamma}_y$ the lift of $\gamma_{\xi}$ starting at $y$.
If $\la g_j \in \psi^{-1}(\xi_0)$ denotes the end point of $\tilde{\gamma}_y$, define
$$
\mu(y) = \la g_j.
$$
Define
$$
K_{U,\la}:=\{(x,y)\in \psi^{-1}(U) \times \psi^{-1}(U) : \ \psi(x) = \psi(y),  \ \  
(\mu(x),\mu(y))-(\la,\la)-1 > 0\}.
$$
and let $\overline{K}_{\la}$ denote the closure of $K_{U,\la}$ in $Y \times Y$. 

\begin{lem} \label{lem3.2}
The divisor $\overline{K}_{\la}$ is an integral symmetric effective correspondence on the 
curve $Y$, 
canonically associated to the triple $(\pi, \mathbb{S}_\omega, \la)$.
\end{lem}

\begin{proof}
For the first assertion it suffices to show that $(\la g, \la)-(\la, \la)-1$ is a non-negative integer for all $g \in W \setminus H$, since
$(\, , )$ is a $W$-invariant scalar product. But $(\la g, \la)-(\la,\la)-1$ is an integer, 
since $\la$ is a weight. Hence it suffices to 
show that $(\la g,\la ) > (\la, \la)$ for every $g \in W \setminus H$. For this note that
$$
(\la , \la) - (\la g,\la) = \frac{1}{2}[(\la g,\la g) -2(\la g,\la ) + (\la ,\la )] = 
\frac{1}{2}(\la g- \la, \la g-\la ),
$$
which implies the assertion, since $( \, , )$ is negative definite and 
$H = \Stab(\la) \subset W$. 

For the last assertion we have to show that the $K_{U,\la}$ does not depend on the 
choice of $\xi_0$ 
the path $\gamma_{\xi}$ connecting $\xi$ and $\xi_0$.
This is a consequence of the $W$-invariance of the form $(\,,)$ and the above mentioned fact that the monodromy map $\rho$
factorizes via the group $W$.
\end{proof}

We call $\overline{K}_{\la}$ the {\it Kanev correspondence} associated to the weight $\la$, 
since it was introduced
by Kanev in \cite{K}.
Considered as a map $Y \ra \Div(Y)$, it is given by
\begin{equation} \label{eq6}
\overline{K}_\la(y)=\sum_{j=1}^d
[(\la g_j, \mu(y))-(\la, \la )-1]\mu^{-1}(\la g_j) + y.
\end{equation}

Note that $y$ is added, since in the sum $y$ appears with
coefficient $-1$, because $(\la g_i,\mu(y))=(\la ,\la )$ if $\mu(y)=\la g_i$.
We need the following description of
$\overline{K}_\la(y)$. According to our construction, for
any $y \in \psi^{-1}(U)$, there is a unique integer $i_y,\; 1 \leq
i_y \leq d$ such that $\mu(y) = \la g_{i_y}$.

\begin{prop}\label{prop3.2}
If $y \in \psi^{-1}(U)$ with $\mu(y) = \la g_{i_y}$, then
$$
\overline{K}_\la (y) = \sum_{j=2}^d [(\la g_j,\la )-(\la ,\la )-1]\mu^{-1}(\la g_jg_{i_y})
$$
\end{prop}

\begin{proof} 
By \eqref{eq6} and using the $W$-invariance of
$(\,,)$, we have
\begin{eqnarray*}
\overline{K}_\la (y) &= & \sum_{\la g_j \neq \mu(y)} [(\la g_j, \mu(y))-(\la,\la)-1]
\mu^{-1}(\la g_j) \\
& = & \sum_{\la g_j \neq \mu(y)}[(\la g_jg_{i_y}^{-1},\mu(y)g_{i_y}^{-1})-(\la ,\la )-1] 
\mu^{-1}(\la g_j) \\
& =&  \sum_{\la g_j \neq \la g_{i_y}} [(\la g_jg_{i_y}^{-1},\la)-(\la,\la)-1] \mu^{-1}(\la g_j)\\
& = & \sum_{j=2}^d [(\la g_j,\la )-(\la,\la)-1] \mu^{-1}(\la g_jg_{i_y}).
\end{eqnarray*}
\end{proof}

\bigskip

\subsection{Relations between $\overline{K}_\la$ and $\overline{S}_\la$}

We want to work out how the Kanev correspondence $\overline{K}_\la$ is related to the Schur correspondence
$\overline{S}_\la$ on the curve $Y$. For this we need a special choice of the set of representatives $g_1, \ldots, g_d$
of the right cosets of the subgroup $H$ of $W$.
In the sequel we choose a set of representatives $\{g_1 =1, \ldots,g_d\}$ such that $\{g_1^{-1}, \ldots,g_d^{-1}\}$
is also a set of representatives of the right cosets of $H$ in $W$. That there is always such a set, is a consequence of the 
marriage theorem of combinatorics (see \cite[Theorem 5.1.7]{hall}).

\begin{prop} \label{prop3.4}
Let $y \in \psi^{-1}(U)$ with $\mu(y) = \la g_{i_y}$. Then
$$\overline{S}_\la (y)=|H|^2\sum_{j=1}^d(\la g_j,\la)\mu^{-1}(\la g_jg_{i_y}).$$
\end{prop}

\begin{proof}
According to Proposition \ref{prop3.1}, 
$\overline{S}_\la (y)=|H|^2\sum_{j=1}^d(\la g_j,\la )\varphi(z^{g_j^{-1}}),$
where $z$ is a point in the fiber $\varphi^{-1}(y) \subset Z$.\\
Suppose first that $\mu(y)=\la $, i.e. $i_y = 1$.
As $W$ acts by right multiplication on the set $\{\la g_1,\ldots, \la g_d\}$, we have 
$\mu(\varphi(z))g = \mu(\varphi(z^g)),$
 and hence
$$
\varphi(z^{g_j^{-1}})=\mu^{-1}(\mu(\varphi(z))g_j^{-1})=\mu^{-1}(\la g_j^{-1})
$$
for all $j$. Therefore we get, using the $W$-invariance of $(\,,)$,
$$
\overline{S}_\la (y) =|H|^2\sum_{j=1}^d(\la g_j,\la )\mu^{-1}(\la g_j^{-1}) 
=|H|^2\sum_{j=1}^d(\la ,\la g_j^{-1})\mu^{-1}(\la g_j^{-1})
=|H|^2\sum_{j=1}^d(\la g_j,\la)\mu^{-1}(\la g_j),
$$
where for the last equation we used the fact that with $\{g_1, \ldots,g_d\}$ also $\{g_1^{-1}, \ldots,g_d^{-1}\}$
is a set of representatives of the right cosets of $H$.\\ 
Finally, if $\mu(y) \neq \la$, we have $\mu(y)g_{i_y}^{-1}= \la$. Hence we can apply the above equation to the bijection
$\widetilde{\mu} = g_{i_y}^{-1} \cdot \mu: \psi^{-1}(\psi(y)) \ra \{\la g_1, \ldots, \la g_d\}$, which gives
$$
\overline{S}_\la (y)=|H|^2\sum_{j=1}^d(\la g_j, \la )\widetilde{\mu}^{-1}(\la g_j)
=|H|^2\sum_{j=1}^d(\la g_j, \la )\mu^{-1}(\la g_jg_{i_y}).
$$
\end{proof}
With the notation and identifications of above we can state the main result of this section

\begin{thm} \label{thm3.5}
The Kanev and Schur correspondence $\overline{K}_\la$ and $\overline{S}_\la$ on the curve
$Y=Z/H$ associated to $\la$ are related as follows
$$
\overline{S}_\la = |H|^2 \left(\overline{K}_\la - \overline{\Delta} + ((\la,\la)+1) \overline{T} \right),
$$
where $\overline{\Delta}$ denotes the diagonal in $Y \times Y$ and $\overline{T} = \psi^*\psi_*$ the trace correpondence of the morphism $\psi$.
\end{thm}

\begin{proof}
It suffices to show that $[\overline{S}_\la + |H|^2 (\overline{\Delta}-\overline{K}_\la)](y)
=[(\la,\la)+1]|H|^2 \overline{T}(y)$ for all
$y \in \psi^{-1}(U)$. But applying Propositions \ref{prop3.2} and \ref{prop3.4} we have
\begin{eqnarray*}
[\overline{S}_\la + |H|^2 (\overline{\Delta}-\overline{K}_\la)](y) &=& 
|H|^2[\sum_{j=1}^d(\la g_j, \la)\mu^{-1}(\la g_jg_{i_y}) +y \\
&& \hspace{2cm}- \sum_{j=2}^d [(\la g_j,\la)-(\la,\la)-1]\mu^{-1}(\la g_jg_{i_y})]\\
&=& |H|^2((\la,\la) +1) \sum_{j=1}^d \mu^{-1}(\la g_jg_{i_y})
\end{eqnarray*}
But the right multiplication with $g_{i_y}$ permutes only the elements $\la g_1, \ldots, \la g_d$, which implies that
$\sum_{j=1}^d \mu^{-1}(\la g_jg_{i_y}) = \psi^{-1}\psi(y) = \overline{T}(y)$. This completes the proof of the theorem.
\end{proof}

\begin{cor}  \label{cor3.6}
 \hspace{2cm}  $\deg \overline{K}_\la = 1 -d((\la,\la )+1)$.
\end{cor}

\begin{proof}
Note first that $\deg S_\la = \sum_{g \in W} (\la g,\la) = (\sum_{g \in W} \la g,\la) = 0$, 
since $\sum_{g \in W} \la g$ is 
$W$-invariant and thus equal to 0 in $V$. This implies $\deg \overline{S}_\la = 0$ and hence Theorem \ref{thm3.5} gives
the assertion, since $\deg \overline{\Delta} = 1$ and $\deg \overline{T} = d$.
\end{proof}

\begin{cor}  \label{cor3.7}
We have the relations 
\begin{enumerate}
\item $\overline{T} \,\overline{S}_\la = \overline{S}_\la \overline{T} = 0 $.
\item $\overline{K}_\la \overline{T} = \overline{T} \, \overline{K}_\la =  \deg \overline{K}_\la \cdot \overline{T}$. 
\end{enumerate}
\end{cor}

\begin{proof}
To obtain (1) we consider $T = \pi^* \pi_*$, the trace correspondence of $\pi$ and note that for all $z \in Z$, 
$$
T S_\la (z) = \sum_{g \in W} (\la g, \la) \,T(z^g) = \sum_{g \in W} (\la g,\la ) \cdot T(z) = \deg S_\la \cdot T(z) = 0
$$
and similarly $S_\la T = 0$. Now $\overline{T} = \frac{1}{|H|} \varphi_* T$ implies the assertion.

To obtain (2) we compute using $\overline{\Delta} \, \overline{T} = \overline{T}$, $\overline{T}^2 = d \overline{T}$
and Theorem \ref{thm3.5}  
$$
\overline{K}_\la \overline{T} = \frac{1}{|H|^2} \overline{S}_\la \overline{T} + \overline{\Delta} \, \overline{T} - ((\la,\la) + 1) 
\overline{T}^2 = \overline{T} - d((\la,\la) + 1)\overline{T} = \deg \overline{K}_\la \cdot  \overline{T},
$$
where the last equation follows from Corollary \ref{cor3.6}.
Similarly Theorem \ref{thm3.5} gives the equation $\overline{K}_\la \overline{T} = \overline{T} \, \overline{K}_\la$.
\end{proof}

For the next theorem we need a lemma.

\begin{lem} \label{lem3.8}
\hspace{2cm} $\overline{S}_\la^2 = \overline{e} \cdot \overline{S}_\la$, with $\overline{e} =   
\frac{|H| \cdot |W|\cdot (\la,\la) }{\dim V_\omega} \in \QQ$.
\end{lem}

\begin{proof}
Schur's orthogonality relations imply that $p_\la = \frac{1}{e} \cdot S_\la$ is an idempotent of the rational group ring 
$\QQ[W]$
or equivalently that $S_\la^2 = e \cdot S_\la$ with $e = \frac{|W|\cdot (\la,\la)}{\dim V_\omega}$, see \cite{LR} 
Proposition 2.3. Note that $e < 0$, since $(\;,\;)$ is negative definite. 

 According to  \cite{F} Proposition 16.1.2 (a), 
$(\varphi \times \varphi)_* S_\la \cdot (\varphi \times \varphi)_* S_\la = 
(\varphi \times \varphi)_* (S_\la \cdot S_\la)$, 
which implies
$$
\begin{array}{ll}
\overline{S}_\la^2 = (\varphi \times \varphi)_* (S_\la^2) & = (\varphi \times \varphi)_* (e \cdot S_\la)\\
& = (\varphi \times \varphi)_* (e \Delta) \cdot (\varphi \times \varphi)_* (S_\la) = |H| e \overline{S}_\la.
\end{array}
$$
\end{proof}

\begin{thm} \label{thm3.9}
The Kanev correspondence $\overline{K}_\la$ satisfies the following cubic equation
$$
(\overline{K}_\la - \overline{\Delta})(\overline{K}_\la + (q_\la - 1) \overline{\Delta})
(\overline{K}_\la - \deg \overline{K}_\la \,\overline{\Delta}) = 0
$$
with $q_\la = - \frac{d \cdot (\la , \la )}{\dim V_\omega} \in \NN$.
\end{thm}

\begin{proof}
Theorem \ref{thm3.5} and Lemma \ref{lem3.8} give 
$$
|H|^4(\overline{K}_\la - \overline{\Delta} + a \overline{T})^2 
= \overline{e} |H|^2(\overline{K}_\la - \overline{\Delta} + a \overline{T})
$$
with $a = (\la ,\la ) +1$. Applying Corollary \ref{cor3.7} and 
$\overline{T}^2 = d \overline{T}$ as well as 
$\overline{\Delta} \, \overline{T} = \overline{T}$, this implies that there is a rational number 
$c$ such that
$$ 
(\overline{K}_\la - \overline{\Delta})^2 - \frac{\overline{e}}{|H|^2}(\overline{K}_\la - \overline{\Delta}) + c \overline{T} = 0.
$$ 
Multiplying this equation by $\overline{K}_\la - \deg \overline{K}_\la \overline{\Delta}$ and using Corollary \ref{cor3.7} again, 
we get
$$
(\overline{K}_\la - \overline{\Delta})(\overline{K}_\la - \overline{\Delta} - 
\frac{\overline{e}}{|H|^2}\overline{\Delta} )
(\overline{K}_\la - \deg \overline{K}_\la \overline{\Delta}) = 0.
$$
This completes the proof of the theorem.
\end{proof}

\bigskip

\subsection{The Prym variety $P_\la$}

We consider the rational correspondence $\overline{s}_\la = \frac{1}{|H|^2}
\overline{S}_\lambda$. Then $\overline{s}_\la$
satisfies the relation $\overline{s}^2_\la =
- q_\la \overline{s}_\la$ and by Theorem 3.5 we have
\begin{equation} \label{endokanevschur}
 \overline{s}_\la = \overline{K}_\la - \overline{\Delta} + ((\la,\la)+1) \overline{T} .
\end{equation}
Hence we see that $\overline{s}_\la$ is an integral correspondence if and only if $(\la,\la)$ is an
integer. This happens e.g. for $E_8$, but not for $E_6$ or $E_7$ -- see the tables in section 4.

\bigskip

The correspondence $\overline{s}_\la$ induces a rational endomorphism also denoted by 
$\overline{s}_\la \in \End_\QQ(\Jac(Y))$. We introduce the following abelian subvariety
of $\Jac(Y)$, which we will also call {\em Prym variety}
$$ P_\lambda := \im ( m \overline{s}_\la ) \subset \Jac(Y),$$
where $m$ is some integer such that $m \overline{s}_\la \in \End(\Jac(Y))$. It is clear 
that $P_\lambda$  does not depend on the integer $m$.

\bigskip

Since $\overline{K}_\la - \overline{\Delta}$ is a correspondence on $Y$ it induces an
endomorphism $v_\la \in \End(\Jac(Y))$. We denote by 
$$ S = \Prym(Y/X) $$ 
the Prym variety of the covering
$\psi : Y \ra X$, see Remark \ref{remarkprym}(ii). Since by Corollary \ref{cor3.7} $v_\la t = t v_\la$,
we obtain that $v_\la(S) \subset S$. Hence $v_\la$ restricts to a symmetric endomorphism
$$ u_\la \in \End(S) \qquad \text{with} \qquad u_\la^2 = q_\lambda u_\la.$$ 

\begin{prop}
We have an equality of abelian subvarieties of $S$
$$ P_\la = \im(u_\la) \subset S.$$
\end{prop}

\begin{proof}
Because of \eqref{endokanevschur} the endomorphisms $\overline{s}_\la$ and $u_\la$ coincide on $S$ 
up to a non-zero integer multiple.
Moreover the restriction of $\overline{s}_\la$ to the subvariety $\psi^* \Jac(X)$ is zero because
$\overline{s}_\la t = 0$ (Corollary \ref{cor3.7}). The equality now follows since 
$\Jac(Y)$ and $\Jac(X) \times S$ are isogenous. 
\end{proof}

\begin{rem}
{\em Note that the abelian subvariety $P_\la \subset S$ only depends on $\la \in V_\omega$
and not on the lattice $\Sbb_\omega \subset V_\omega$. Moreover we recall from \cite{M}
Proposition 4.3 (1) that $P_\la$ and $P_{\la'}$ are isogenous for any $\la,\la' \in 
V_\omega$.}
\end{rem}

\begin{rem} 
{\em If $\psi: Y \ra X$ is \'etale, then the Kanev correspondence $\overline{K}_\la$ is 
fixed-point-free. On the other hand we notice that the induced endomorphism $v_\la$ does not satisfy a quadratic equation --- only its restriction $u_\la$ to $S$ --- hence we are not
in a situation where Kanev's criterion (\cite{K2} or \cite{BL}  Theorem 12.9.1) applies.}
\end{rem}

\section{Examples}

As our main examples let us work out the invariants of the last section in the case of a Weyl group $W$ of ADE-type. Let 
$\pi: Z \ra X$ be a Galois covering of smooth projective curves with Galois group $W$. 
As weights $\lambda$ we choose the 
minuscule weights $\varpi_i$ (with the notation of \cite{Bou}) except in the case of type $E_8$, where we consider the 
quasi-minuscule weight $\lambda = \varpi_8$. As above let $V_{\lambda}$ denote the corresponding absolutely irreducible 
$\QQ$-representation, $H$ the stabilizer subgroup of $\lambda$ in $W$ and $Y = Z/H$. If $(\;, \;)$ denotes the 
uniquely determined $W$-invariant negative definite symmetric form on $V_{\lambda}$ as defined
in section 3.2, the following tables give 
the values for the numbers $(\lambda,\lambda),\; d = \deg(Y/X), \;\dim V_{\lambda},\; q_{\lambda}$ and the degree 
of the Kanev-correspondence $\overline{K}_{\lambda}$. For the definition of the Dynkin index $d_\la$ see 
\cite{LS} and also Proposition \ref{propls}. The case of a Weyl group of type $E$ is closely related to a del 
Pezzo surface (see \cite{K}). 
The last line of the first table gives the degree of the corresponding surface.\\

\begin{center}
\begin{tabular}{|c||c|c|c|c|c|}
 Weyl group of type & $E_4 = A_4$ & $E_5 = D_5$ & \quad $E_6$ \quad \, & \quad $E_7$ \quad \, & \quad $E_8$ \quad \,\\ \hline \hline
weight $\lambda $ & $\varpi_2$ & $\varpi_4, \varpi_5$ & $\varpi_1, \varpi_6$ & $\varpi_7$ & 
$\varpi_8$\\ \hline 
$(\lambda,\lambda)$ & $- \frac{6}{5}$ & $- \frac{5}{4}$ & $- \frac{4}{3}$ & $- \frac{3}{2}$ & $-2$ \\ \hline
$d = \deg (Y/X)$ & 10 & 16 & 27 & 56 & 240 \\ \hline 
$\dim V_{\lambda}$ & 10 & 16 & 27 & 56 & 248 \\ \hline
$ q_{\lambda} = d_{\lambda}$ & 3 & 4 & 6 & 12 & 60 \\ \hline 
$\deg \overline{K}_{\lambda}$ & 3 & 5 & 10 & 29 & 241 \\ \hline 
del Pezzo of degree & 5 &4&3&2&1\\ \hline

\end{tabular}
\end{center}

\vspace{0.4cm}

\begin{center}
\begin{tabular}{|c||c|c|}
 Weyl group of type & $A_n$ & $D_n$\\ \hline \hline
 weight $\lambda$ & $\varpi_i, \, 1 \leq i \leq n$ & $\varpi_{n-1}, \varpi_n$ \\ \hline
$(\lambda,\lambda)$ & $- \frac{i(n+1-i)}{n+1}$ & $- \frac{n}{4}$ \\ \hline 
$d = \deg (Y/X)$ & ${n+1 \choose i}$ & $ 2^{n-1}$ \\ \hline 
$\dim V_{\lambda}$ & ${n+1 \choose i}$ & $2^{n-1}$ \\ \hline
$ q_{\lambda} = d_{\lambda}$ & ${n-1 \choose i-1}$ & $2^{n-3}$ \\ \hline 
$\deg \overline{K}_{\lambda}$ & ${n-1 \choose i-1}n - {n+1 \choose i } +1$ & $2^{n-3}(n - 4) + 1$        \\ \hline

\end{tabular}
\end{center}

\bigskip

\begin{proof}
The symmetric form $(\;,\;)$ is given by the negative of the Cartan matrix. This gives the value of $(\lambda,\lambda)$
using the explicit form of $\lambda$ as outlined in \cite{Bou}. The length of the orbit of $\lambda$ under the action of $W$ 
gives the degree of $Y/X$. The dimension of $V_{\lambda}$ follows from the fact that $\lambda$ is a minuscule weight, 
respectively quasi-minuscule weight in the case of $E_8$. The values of $q_{\lambda}$ and $\deg \overline{K}_{\lambda}$ 
are computed by Theorem 3.9 and Corollary 3.6. For the degree of the del Pezzo surface see \cite{K}. For the Dynkin indices see \cite{LS} Proposition 2.6.
\end{proof}

\begin{rem}
{\em  We note that in the examples of the tables $q_\la = d_\la$. This equality of the two
integers is a coincidence, as e.g. for $(G_2, \varpi_1)$ : $d_\la = 2$, $q_\la = 6$ or for
$(F_4, \varpi_4)$ : $d_\la = 6$, $q_\la = 12$ see \cite{LS} Proposition 2.6 and \cite{K} 
page 176.
}
\end{rem}

\bigskip

\section{The Donagi-Prym variety}

\subsection{Definitions and properties}

Following \cite{D} we introduce the compact commutative algebraic group 
$$ \Prym(\pi, \Sbb_\omega) := \Hom_W(\Sbb_\omega, \Pic(Z)), $$
which we call the Donagi-Prym variety associated to the pair $(\pi, \Sbb_\omega)$,
where $\pi: Z \ra X$ is a Galois covering with group $W$ and $\Sbb_\omega$ is a 
right-$\ZZ[W]$-module.  
The elements of $\Prym(\pi, \Sbb_\omega)$ are $W$-equivariant 
homomorphisms $\phi: \Sbb_\omega \ra \Pic(Z)$; for an
alternate description see Lemma \ref{bijdp}. Its connected component $\Prym(\pi,\Sbb_\omega)_0$
containing $0$ is an abelian variety. For any element $\lambda \in 
\Sbb_\omega$ we consider the corresponding evaluation map 
$$ ev_\lambda : \Prym(\pi, \Sbb_\omega) \lra \Pic(Z),  \qquad \phi \mapsto \phi(\lambda). $$

\begin{prop}
We have the equality
$$ \Hom_W (\Sbb_\omega, \Pic(Z) ) = \Hom_W (\Sbb_\omega, \Jac(Z) ).$$
In particular all connected components of the image
$ev_\lambda ( \Prym(\pi, \Sbb_\omega) )$ are contained in $\Jac(Z)$.
\end{prop}

\begin{proof}
Consider a homomorphism $\phi \in \Hom_W (\Sbb_\omega, \Pic(Z) )$. The composite map $\deg \circ
\phi$ is a $W$-invariant homomorphism from $\Sbb_\omega$ to $\ZZ$. Since the dual lattice
$\Hom_\ZZ(\Sbb_\omega, \ZZ)$ is an irreducible $W$-module, we obtain that   $\deg \circ
\phi = 0$. 
\end{proof}

\bigskip

We denote by 
$$\Gamma_\lambda = \Sbb_\omega/ (\lambda \cdot \ZZ[W]).$$ 
This is a finite abelian group since 
$V_\omega = \Sbb_\omega \otimes_\ZZ \QQ$ is an irreducible $W$-module.

\begin{prop} \label{kerevla}
The kernel of the evaluation map equals 
$$ \ker ev_\lambda = \Hom_W (\Gamma_\lambda, \Jac(Z) ). $$
\end{prop}

\begin{proof}
Let $U$ denote the image of the natural map $\ZZ[W] \stackrel{a}{\lra}  U  \stackrel{i}{\lra} 
\Sbb_\omega$ 
defined by $a(g)= \lambda g$ for $g \in W$. Then we have
the exact sequence of right $\ZZ[W]$-modules
$$ 0 \lra U \stackrel{i}{\lra} \Sbb_\omega \lra \Gamma_\lambda \lra 0.$$
Now we apply the left-exact contravariant functor $\Hom_W(\cdot, \Jac(Z))$ and we obtain
the exact sequence
$$ 0 \lra \Hom_W(\Gamma_\lambda, \Jac(Z)) \lra \Hom_W(\Sbb_\omega, \Jac(Z)) 
\stackrel{\check{i}}{\lra} \Hom_W(U, \Jac(Z)).$$
We then apply the same functor to the surjective $\ZZ[W]$-morphism $a: \ZZ[W] \lra U$ and we obtain
an injective morphism
$$ \check{a}: \Hom_W(U, \Jac(Z)) \lra \Hom_W(\ZZ[W], \Jac(Z)) = \Jac(Z).$$
Note that the composite morphism $\check{a} \circ \check{i}$ equals the evaluation map 
$ev_\lambda$ and
that $ \ker ev_\lambda = \ker \check{i}$ because $\check{a}$ is injective. 
\end{proof}

\bigskip

\subsection{The image of the evaluation map $ev_\la$}

It is clear that 
$$ev_\la (\Prym(\pi, \Sbb_\omega)) \subset \Jac^H(Z),$$
where $H = \mathrm{Stab}(\la)$ and $\Jac^H(Z)$ denotes the subvariety of $\Jac(Z)$ 
parametrizing $H$-invariant line bundles over $Z$. We recall the exact sequence of
abelian groups (\cite{Do} Proposition 2.2)
$$
0 \lra H^* \lra \Pic(H;Z) \map{\iota} \Pic^H(Z) \lra H^2(H,\CC^*) \lra 0,
$$
where $\Pic(H;Z)$ denotes the group of $H$-linearized line bundles over $Z$ (see 
\cite{Do} section 1) and $H^* = \Hom(H,\CC^*)$ denotes the group of
characters of $H$. The map $\iota$ is the map which forgets the $H$-linearization. 

\begin{prop}
We have the exact sequence 
\begin{equation} \label{esetale}
0 \lra H^* \lra \Jac(H;Z) \map{\iota} \Jac^H(Z) \lra H^2(H,\CC^*) \lra 0.
\end{equation}
\end{prop}

\begin{proof}
We denote by $d$ the degree of the cover $\psi: Y \rightarrow X$. Consider the
following diagram, in which the first line is obtained as the kernel of the degree
morphism.
$$
\begin{array}{ccccccc}
 & 0           & &  0        &  & 0 & \\
 & \downarrow &  & \downarrow & & \downarrow &  \\
0 \ra & \Jac(Y)/H^* & \ra & \Jac^H(Z) & \stackrel{c}{\ra} & H^2(H, \CC^*) & \\
 & \downarrow & & \downarrow & & \downarrow \scriptstyle{\cong} & \\
0 \ra & \Pic(Y)/H^* & \ra & \Pic^H(Z) & \stackrel{c}{\ra} & H^2(H, \CC^*) &  \ra 0\\
 & \downarrow \scriptstyle{\deg} & & \downarrow \scriptstyle{\deg} & & \downarrow 
\scriptstyle{0} & \\
0 \ra & \ZZ & \stackrel{\cdot d}{\ra} & \ZZ & \ra & \ZZ/ d \ZZ & \ra 0 \\
& \downarrow & & \downarrow & & \downarrow \scriptstyle{\cong} & \\
& 0           & \ra & \ZZ/d \ZZ & \stackrel{\cong}{\ra} & \ZZ/ d \ZZ & \ra 0 
\end{array}
$$
The snake lemma now implies that the natural homomorphism $c : \Jac^H(Z) \ra H^2(H, \CC^*)$ is surjective.
\end{proof}

In section 6.5
we will show that for any $\la \in \Sbb_\omega$ the evaluation map
$ev_\la$ lifts  to $\Jac(H; Z)$. 
Finally we also recall the exact sequence 
\begin{equation} \label{ramZY}
0 \lra \Jac(Y) \lra \Jac(H; Z) \lra \bigoplus_{x \in \mathrm{Br}(\varphi)}
\ZZ/ e_x \ZZ \lra 0,
\end{equation}
where $x$ varies over the branch divisor $\mathrm{Br}(\varphi)$ of the 
Galois covering $\varphi : Z \ra Y = Z/H$ and $e_x$ denotes the 
ramification index of $x$.

\bigskip

\subsection{The group $\Gamma_\la$}

We now compute the finite group $\Gamma_\la$ in some special cases, which will
be worked out in detail in section 8.  Given a simple
Lie algebra $\gfr$ we take for $\Sbb_\omega$ the weight lattice of $\gfr$ and
denote by $\Lambda_R \subset \Sbb_\omega$ the root lattice.
We use the notation of \cite{Bou}. In particular $\alpha_i$ denote the simple roots
and $\varpi_i$ denote the fundamental weights of $\gfr$.

\begin{lem} \label{groupgamma}
In the following cases the group $\Gamma_\lambda$ is trivial.
\begin{enumerate}
\item $W = W(E_8)$, $\lambda = \varpi_8$.
\item $W = W(E_7)$, $\lambda = \varpi_7$.
\item $W = W(E_6)$, $\lambda = \varpi_1$ or $\lambda = \varpi_6$.
\item $W = W(D_n)$ with $n$ odd, $\lambda = \varpi_{n-1}$ or $\lambda = \varpi_n$.
\item $W = W(A_n)$, $\lambda =  \varpi_i$ with $i$ coprime to $n+1$.
\end{enumerate}
\end{lem}

\begin{proof}
In each case we have to show that $\lambda \cdot \ZZ[W] = \Sbb_\omega$. Let $\alpha \in
\Lambda_R$ be a root and let $s_\alpha$ be its associated reflection. We have
$$ s_\alpha(\varpi) = \varpi - (\varpi | \alpha) \alpha \qquad \forall \varpi \in 
\Sbb_\omega.$$
Note that for Lie algebras of type ADE all roots $\alpha$ have same length $(\alpha |\alpha) 
= 2$. Here $( \cdot | \cdot )$ denotes the Cartan-Killing form on $\Lambda_R$.

\begin{enumerate}
\item Clearly the elements $\varpi_8$ and $s_{\alpha_8}(\varpi_8)$ are in the lattice
$\lambda \cdot \ZZ[W]$. Moreover $s_{\alpha_8}(\varpi_8) = \varpi_8 - \alpha_8$, hence 
$\alpha_8 \in \lambda \cdot \ZZ[W]$. The Weyl group $W$ acts transitively on the 
roots, which implies that $\Lambda_R \subset \lambda \cdot \ZZ[W]$. For $\gfr$ of type $E_8$, we have $\Lambda_R = \Sbb_\omega$ and we are done.

\item As before $\varpi_7, s_{\alpha_7}(\varpi_7) \in \lambda \cdot \ZZ[W]$, hence
$\alpha_7 \in  \lambda \cdot \ZZ[W]$. Since $W$ acts transitively on the roots, we have
$\Lambda_R \subset \lambda \cdot \ZZ[W]$. Here $\varpi_7 \in \la \cdot \ZZ[W]$, but 
$\varpi_7 \notin \Lambda_R$. Since $[\Sbb_\omega : \Lambda_R] = 2$, this implies
the assertion.

\item The computations are similar to the previous case using $(\alpha_1 | \varpi_1) = 
(\alpha_6 | \varpi_6) = 1$ and $[\Sbb_\omega : \Lambda_R] = 3$.

\item Since $W$ acts transitively on the roots, we obtain as before that $\Lambda_R \subset
\lambda \cdot \ZZ[W]$. Since $n$ is odd, we have that $\Sbb_\omega/\Lambda_R \cong \ZZ/ 4 \ZZ$ is
generated by the class of $\varpi_n$ or $\varpi_{n-1}$.

\item As before we obtain that $\Lambda_R \subset \lambda \cdot \ZZ[W]$. The class of the 
fundamental weight $\varpi_i$ in the quotient $\Sbb_\omega/ \Lambda_R \cong \ZZ / (n+1) \ZZ$ equals
the class of $i \in \ZZ / (n+1) \ZZ$. The assertion then follows. 
\end{enumerate}
\end{proof}

\bigskip

\section{Mumford groups}

\subsection{Twisted Mumford groups}

Let $T$ be a torus and suppose we are given a left action of the group $W$ on $T$
$$ \sigma : W \lra \Aut(T).$$
For any $g \in W$  we also denote the automorphism of $T$ by $g$. This left action $\sigma$
induces a right action of $W$ on the group of characters $\Hom(T, \CC^*)$. We suppose that
the representation of $W$ on $\Hom(T, \CC^*) \otimes_\ZZ \QQ$ is irreducible, hence 
$\Hom(T, \CC^*)$ is
of the form $\Sbb_\omega$ for some $\omega \in \widehat{W}$. Conversely given a lattice $\Sbb_\omega$
there always exists a torus $T$ with a left $W$-action such that
$$ \Sbb_\omega := \Hom(T, \CC^*)$$
as $W$-modules. We recall that $W$ acts from the right on the curve $Z$. 

\bigskip

Let $E_T$ be a principal $T$-bundle over $Z$. The group $W$ then acts from the left
on the set of principal $T$-bundles in two different ways: $g \in W$ sends
$E_T$ to
\begin{enumerate}
\item $g^*E_T$, the pull-back of $E_T$ under the automorphism $g$ of $Z$.
\item $E_T  \times^T_g T$, the $T$-bundle obtained from $E_T$ by extension of structure group from $T$ to $T$, where $t \in T$ acts on $T$ by left-multiplication with $g(t) \in T$.
\end{enumerate}

Note that both actions commute. In this paper we will denote the combined left $W$-action on the set of principal $T$-bundles by
$$g. E_T = g^* E_T \times^T_g T, \qquad \text{for} \ g \in W.$$

Since the action of $W$ on $Z$ is a right action, we have $(g_1 g_2)^* E_T = g_1^* g_2^* E_T$.
This implies that the action $(g, E_T) \mapsto g . E_T$ is a left action.

\bigskip

We say that a $T$-bundle $E_T$ is $W$-{\em invariant} if $E_T \simeq g.E_T$ for all $g \in W$.
Here the symbol $\simeq$ denotes a $T$-equivariant isomorphism. Let $\underline{T}$ denote the
sheaf of abelian groups on $Z$ defined by $\underline{T}(U) = \Mor(U, T)$ for any open subset 
$U \subset Z$. The cohomology group $H^1(Z,\underline{T})$ parametrizes 
isomorphism classes of $T$-bundles over $Z$.  We denote by  $H^1(Z,\underline{T})^W$
the $W$-invariant subset, where $W$ acts both on the curve $Z$ and on the torus $T$.

\begin{lem} \label{bijdp}
There are canonical bijections
\begin{eqnarray*}
H^1(Z,\underline{T})^W & = & \{ W \text{-invariant} \;  T\mbox{-bundles} \; E_T \; \mbox{on} \; Z \; \} \\
 & = &  \Prym(\pi, \Sbb_\omega).
\end{eqnarray*}
\end{lem}
\begin{proof}
The first equality follows from the definition of $H^1(Z,\underline{T})^W$.  For the second equality  note that the map 
$$
\Phi: \{ T\mbox{-bundles} \; E_T \; \mbox{on} \; Z \} \lra \Hom( \Sbb_\omega , \Pic(Z))
$$
$$
\Phi(E_T) =  (\phi: \Sbb_\omega \ra \Pic(Z), \; \phi(\lambda) = E_T \times^T_{\lambda} \CC^*)
$$
is bijective. The inverse map is given by mapping $\phi \in \Hom( \Sbb_\omega , \Pic(Z))$ 
to $E_T = L_1 \oplus \cdots \oplus L_n$, where 
$L_i = \phi(\lambda_i)$, $\lambda_1, \ldots, \lambda_n$ being a $\ZZ$-basis of $\Sbb_\omega$.
One can easily show that the inverse map does not depend on the choice of the $\ZZ$-basis.

We have to show that $\phi = \Phi(E_T)$ is $W$-equivariant if and only if $E_T \simeq  g.E_T$ 
for all $g \in W$. Since 
$W$ acts from the right on $\Sbb_\omega$ and $Z$, the map $\phi$ is $W$-equivariant if and only if 
\begin{equation} \label{eqn}
\phi(\lambda g) = (g^{-1})^*\phi(\lambda)
\end{equation}
for all $\lambda \in \Sbb_\omega$ and $g \in W$. But 
$$ 
\phi(\lambda g) = E_T \times^T_{\lambda g}\CC^* = (E_T \times^T_g T) \times^T_{\lambda} \CC^*
$$
and 
$$ 
(g^{-1})^*\phi(\lambda) = (g^{-1})^*(E_T \times^T_{\lambda} \CC^*).
$$
Hence (\ref{eqn}) is valid for all $\lambda \in \Sbb_\omega$ if and only if 
$(g^{-1})^*E_T \simeq E_T \times^T_g T$, 
which is the case if and only if $E_T \simeq g^*E_T \times^T_{g} T = g. E_T$.
\end{proof}

\begin{defi}
{\em Given a $W$-invariant $T$-bundle $E_T$ over $Z$ we define the $\sigma$-{\em twisted Mumford group of} $E_T$ by
$$ \cG^\sigma(E_T) := \{ (\gamma, g) \ | \ g \in W, \ \gamma: E_T \lra g. E_T \},$$
where $\gamma$ is a $T$-equivariant isomorphism. The composition law in 
$\cG^\sigma(E_T)$ is given by
$$ (\gamma, g) . (\mu,g') = ((g. \mu) \circ \gamma, gg'), $$
where $g. \mu : g. E_T \lra g. (g'. E_T) = (gg'). E_T$ denotes the isomorphism obtained by applying $g.$ to $\mu$.}
\end{defi}

\bigskip

Since the automorphism group $\Aut_T(E_T)$ equals $T$, we obtain that the group $\cG^\sigma(E_T)$ fits into the exact sequence
\begin{equation} \label{extmg}
 1 \lra T \lra  \cG^\sigma(E_T) \lra W \lra 1,
\end{equation}
where the second map is the map which forgets the isomorphism $\gamma$. We note 
that the induced action by conjugation of $W$ 
on the normal subgroup $T$ coincides with $\sigma$. To simplify notation 
we denote the $\sigma$-twisted Mumford group $\cG^\sigma(E_T)$ by $N$.

\begin{prop} \label{actionN}
Let $E_T$ be a $W$-invariant $T$-bundle over $Z$ with projection $p: E_T \ra Z$. 
Then the natural $T$-action on $E_T$
extends canonically to a $N$-action on the variety $E_T$ preserving the fibers of 
$\pi \circ p: E_T \ra X$, i.e. there exists a morphism over $X$ defining a group action
$$ \mu : E_T \times N \lra E_T, $$
which extends right-multiplication of $T$ on $E_T$.
\end{prop}

\begin{proof}
Given $n = (\gamma, g) \in N$ we first define 
$$\nu_n := \beta_g \circ \gamma,$$ 
where $\beta_g$ is the pull-back of $g : Z \ra Z$ via the projection map
$p: E_T \times^T_g T \ra Z$. This is
summarized by the diagram
$$
\begin{array}{ccc}
E_T &  &  \\
  & & \\
 \downarrow {\scriptstyle \gamma} & \searrow {\scriptstyle \nu_n}  &    \\
   & & \\
g^* (E_T \times^T_g T) & \map{\beta_g} & E_T \times^T_g T \\
 & & \\
\downarrow & & \downarrow {\scriptstyle p} \\
 & & \\
Z & \map{g} & Z.
\end{array}
$$
In order to define the morphism $\mu$, it suffices to define
for any $n=(\gamma, g) \in N$ the right-multiplication with $n$ on
$E_T$
\begin{equation} \label{multn}
 \mu_n : E_T \lra E_T \qquad \text{satisfying} \qquad \mu_{nn'} = \mu_{n'} \circ
\mu_n \ \  \forall n,n' \in N.
\end{equation}
The set of morphisms $\{ \mu_n \}_{n \in N}$ is related to
the set of ``twisted" morphisms $\{ \nu_n \}_{n \in N}$ by composition

%Since $ntn^{-1} = g(t)$, we also have $n^{-1} t n = g^{-1}(t)$. Hence for all 
%$e \in E_T$, $n= (\gamma, g) \in N$ and $t \in T$ we can write
%$$ etn = en g^{-1}(t).$$
%This equality implies that the morphism
%$$ \nu_n : E_T \lra E_T \times^T_g T, \qquad e \mapsto \overline{(en,1)}$$
%is a $T$-equivariant isomorphism. Here $ \overline{(en,1)}$ denotes the class
%of $(en,1) \in E_T \times T$ in the quotient $E_T \times^T_g T := E_T \times T / T$.
%Conversely given the set $\{  \nu_n \}_{n \in N}$ we recover the set 
%$\{ \mu_n \}_{n \in N}$ by composing

$$ \mu_n : E_T \map{\nu_n} E_T \times^T_g T \map{\alpha_g} E_T, \qquad \text{with}
\ \alpha_g \left[ \overline{(e,t)}  \right] = e g^{-1}(t). $$

Here $ \overline{(e,t)}$ denotes the class
of $(e,t) \in E_T \times T$ in the quotient $E_T \times^T_g T := E_T \times T / T$. Note
that we use the equality $etn = en g^{-1}(t)$ for $e \in E_T$, $n= (\gamma, g) \in N$ 
and $t \in T$.

\bigskip

It is now straightforward to check that the set $\{ \mu_n \}_{n \in N}$ satisfies
\eqref{multn}.
\end{proof}

\bigskip

\subsection{The relative case}

In order to construct the abelianization map (see section 7.3) we need to define the
twisted Mumford group for families of $W$-invariant $T$-bundles.

\bigskip

Let $S$ be a scheme and let $\cE_T$ be a family of $T$-bundles over $Z$ parametrized by
$S$, i.e. a $T$-bundle over the product $Z \times S$. Let $\pi_S: Z \times S \ra S$ denote
projection onto the second factor. We say that $\cE_T$ is $W$-{\em invariant} if for any 
$g \in W$ and any closed point $s \in S$ there exists an isomorphism
\begin{equation} \label{famET}
 {\cE_T}_{| Z \times \{ s \} } \simeq g. {\cE_T}_{| Z \times \{ s \} }.
\end{equation}

The following lemma is the analogue of the see-saw theorem for $T$-bundles (see \cite{Mu}).

\begin{lem}
Let $\cA_T$ be a family of $T$-bundles over $Z$ parametrized by $S$. 
Suppose that for any closed point $s \in S$ the $T$-bundle
$ {\cA_T}_{| Z \times \{ s \} }$ is trivial. Then there exists a $T$-bundle $M$ 
over $S$ such that $\cA_T \simeq \pi_S^* M$.
\end{lem}

Applying the preceding lemma to \eqref{famET} we obtain that for any $g \in W$
there exists a $T$-bundle $M_g$ over $S$ such that
$$ \cE_T \simeq g. \cE_T \otimes \pi_S^* M_g .$$
Moreover the map $W \ra H^1(S,\underline{T})$, $g \mapsto M_g$ is easily seen to 
be a group homomorphism.

\bigskip

We now define the $\sigma$-twisted Mumford group of the family $\cE_T$ by
$$ \cG^\sigma(\cE_T) := \{ (\gamma, g) \ | \ g \in W, \ \gamma: 
\cE_T \lra g. \cE_T  \otimes \pi_S^* M_g \},$$
where $\gamma$ is a $T$-equivariant  isomorphism over $Z \times S$. We easily see that
the relative Mumford group fits in the exact sequence 
$$ 1 \lra \Mor(S,T) \lra \cG^\sigma(\cE_T)  \lra W \lra 1.$$
If $S$ is complete and connected, then $\Mor(S,T) = T$ and for any closed point $s \in S$ the
restriction map
\begin{equation} \label{resmumford}
 \mathrm{res}_s: \  \cG^\sigma(\cE_T)  \map{\sim} \cG^\sigma(E_T),
\qquad (\gamma,g) \mapsto ( \mathrm{res}_s(\gamma) , g)
\end{equation}
is an isomorphism. Here $E_T :=  {\cE_T}_{| Z \times \{ s \} }$.

\bigskip

\subsection{Non-twisted Mumford groups}

Let $G$ be an affine algebraic group.

\begin{defi}
{\em Given a principal $G$-bundle $E_G$ over $Z$ we define the (non-twisted) Mumford group
of $E_G$ by
$$ \cG^1(E_G) := \{ (\gamma, g) \ | \ g \in W, \ \gamma: E_G \lra g^* E_G \},$$
where $\gamma$ is a $G$-equivariant isomorphism. The composition law in 
$\cG^1(E_G)$ is given by
$$ (\gamma, g) . (\mu,g') = (g^*(\mu) \circ \gamma, gg'). $$
}
\end{defi}

As before, the group $\cG^1(E_G)$ fits into the exact sequence 
\begin{equation} \label{exmg}
1 \lra \Aut(E_G) \lra \cG^1(E_G) \lra W \lra 1.
\end{equation}
The induced action by conjugation of $W$ on $\Aut(E_G)$ is trivial.

\begin{prop} \label{splitEN}
Let $E_T$ be a $W$-invariant $T$-bundle over $Z$ and denote $N = \cG^\sigma(E_T)$. Let
$E_N := E_T \times^T N$ be the principal $N$-bundle obtained by extending the structure 
group from $T$ to $N$. Then the exact sequence \eqref{exmg} for $E_N$
$$ 1 \lra \Aut(E_N) \lra \cG^1(E_N) \lra W \lra 1 $$
splits canonically.
\end{prop}

\begin{proof}
By Proposition \ref{actionN} the variety $E_T$ admits a canonical $N$-action. Given
$e \in E_T$ and $n = (\gamma, g) \in N$ we denote their product by $en$. We consider
the morphism
$$ \Phi_n : E_T \times^T N \lra E_T \times^T N, \qquad \Phi_n \left[ 
\overline{(e,f)}   \right] = \overline{(en,n^{-1}f)} $$
with $e \in E_T$ and $f \in N$. Then one easily checks that $\Phi_n$ is
well-defined and that it makes the following diagram commute
$$
\begin{CD}
E_T \times^T N @>\Phi_n>> E_T \times^T N \\
@VVpV    @VVpV \\
Z @>g>> Z.
\end{CD}
$$
Hence $\Phi_n$ is a lift of the automorphism $g$ to the principal $N$-bundle 
$E_T \times^T N$ and we may view $\Phi_n$ as an element of $\cG^1(E_N)$.
Moreover $\Phi_{nn'} = \Phi_{n'} \circ \Phi_n$ and $\Phi_t = \mathrm{id}$ for all
$n,n' \in N$ and all $t \in T$. Hence the morphism $\Phi_n$ only depends on $g \in W$. This
gives the canonical splitting.
\end{proof}

\subsection{Mumford groups and descent}  

We recall the relation between the non-twisted Mumford group $\cG^1(E_G)$ and descent.

\begin{defi}
Let $E_G$ be a $G$-bundle over $Z$. A $W$-linearization of the $G$-bundle $E_G$ is a
splitting of the exact sequence \eqref{exmg}
$$ \cG^1(E_G) \map{\curvearrowleft} W. $$
\end{defi}

Descent theory gives in the usual way (see \cite{Gr2})

\begin{prop} \label{descentEG}
Let $E_G$ be a $W$-linearized $G$-bundle over $Z$. There exists a $G$-bundle $F_G$ over
$X$ such that 
$$ E_G \simeq \pi^* F_G $$
if and only if for any ramification point $z \in Z$ of the covering $\pi: Z \ra X$
the stabilizer subgroup $W_z \subset W$ of $z$ acts trivially on the fiber $(E_G)_z$.
\end{prop}

\begin{cor}
If the covering $\pi: Z \ra X $ is \'etale, then given $E_G$ over $Z$ there is a
bijection between the set of $W$-linearizations on $E_G$ and the set of $G$-bundles $F_G$ over
$X$ such that $E_G \simeq \pi^* F_G$.
\end{cor}

\begin{rem}
{\em  We can rephrase Proposition \ref{splitEN} by saying that the $N$-bundle $E_N :=
E_T\times^T N$
admits a canonical $W$-linearization. The conditions for descent of $E_N$ to $X$  given
in Proposition \ref{descentEG} then translate into certain conditions on the $W$-invariant
$T$-bundle $E_T$ --- for the case of cameral coverings, see \cite{DG}.
}
\end{rem}

\subsection{The evaluation map $ev_\la$ lifts}

\begin{prop} \label{evlift}
For any $\la \in \Sbb_\omega$ the evaluation map 
$ev_\la$ lifts to $\Jac(H; Z)$; i.e. there is a commutative diagram
$$
\begin{array}{ccc}
   &   & \Jac(H;Z) \\
   &   &   \\
   & {\scriptstyle \widetilde{ev}_\la}  \nearrow   &  \downarrow \iota \\ 
   &   &   \\
\Prym(\pi, \Sbb_\omega) & \map{ev_\la} & \Jac^H(Z).
\end{array}
$$
\end{prop}

\begin{proof}
We have to show that the line bundle $L = ev_\la(E_T) = E_T \times^T_\la \CC^*$ admits a
canonical $H$-linearization. We denote by $\la = \la_1, \ldots , \la_d \in \Sbb_\omega$ 
the $W$-orbit of the element $\la \in \Sbb_\omega$, with $d = [W: H]$ and by 
$L_i := E_T \times^T_{\la_i} \CC^* \in \Jac(Z)$ the line bundle obtained through $\la_i$
for $ i= 1, \ldots,d$. As before we denote $N = \cG^\sigma(E_T)$ and we consider
$V = \Ind_T^N(\CC_\la)$ the induced representation of $N$ from the representation 
$T \ra GL(\CC_\la)$ given by $\la \in \Hom(T,\CC^*)$. The restricted representation
decomposes as a $T$-module (\cite{Se} Proposition 15)
$$\Res_T (V) = \left[ \bigoplus_{i=1}^d \CC_{\la_i} \right]^{\oplus |H|} $$
which implies that the vector bundle $E_N \times^N V$ decomposes as 
$\left[ \bigoplus_{i=1}^d L_i \right]^{\oplus |H|}$. Moreover by Proposition 
\ref{splitEN} $E_N$ admits a canonical $W$-linearization, hence the decomposable
bundle $E_N \times^N V$ also does. Since the subgroup $H \subset W$ 
preserves the direct summand $L = L_1$, we are done.
\end{proof}

\bigskip

In section 7.3 we will need the next lemma. First we introduce some notation. 
Given $\la \in \Sbb_\omega$ we denote by $\la = \la_1, \ldots, \la_d \in \Sbb_\omega$
the $W$-orbit of $\la$ and by $L_i := E_T \times^T_{\la_i} \CC^*$ the
associated line bundles for $i= 1, \ldots , d$. We recall
(see \eqref{ramZY})  that if $\pi : Z \ra X $ is \'etale, then $\Jac(Y) = 
\Jac(H; Z)$  and we denote by $M \in \Jac(Y)$ the line bundle defined by the 
relation 
$$ \varphi^* M =  L_1 = ev_\la (E_T). $$

\begin{lem} \label{dirim}
The direct sum $\bigoplus_{i=1}^d L_i$ admits a canonical $W$-linearization and, if 
$\pi: Z \ra X$ is \'etale, then the $W$-linearized vector bundle $\bigoplus_{i=1}^d L_i$
equals $\pi^*(\psi_* M)$.
\end{lem} 

\begin{proof}
We have already seen in the proof of Proposition \ref{evlift} that $\bigoplus_{i=1}^d L_i$
admits a canonical $W$-linearization. If $\pi: Z \ra X$ is \'etale, then there exists 
a vector bundle $A$ over $X$ such that $\pi^* A  = \bigoplus_{i=1}^d L_i$ as
$W$-linearized vector bundles. By adjunction for the morphism $\psi: Y \ra X$ we have
$\Hom(A, \psi_* M) = \Hom( \psi^* A, M)$ and since $\pi$ is \'etale, we have
$$\Hom( \psi^* A, M) = \Hom( \pi^* A , \varphi^* M) ^H = \Hom ( \bigoplus_{i=1}^d L_i ,
L_1)^H,$$
where the exponent $H$ denotes $H$-equivariant homomorphisms. The last space 
contains the projection onto the first factor, which is $H$-equivariant. This 
gives a nonzero homomorphism $\alpha : A \ra \psi_* M$.

\bigskip

Similarly one can show using the isomorphism $\mathcal{H}om(\psi_* M,\cO_X) = \psi_* (M^{-1})$ that
\\ 
$ \Hom(\psi_* M, A) = \Hom( \bigoplus_{i=1}^d L^{-1}_i , L^{-1}_1)^H$. As before
projection onto the second factor gives a nonzero homomorphism $\beta: \psi_*M \ra A$
and one easily checks that $\alpha \circ \beta = \mathrm{id}_{\psi_*M}$. Hence
$\alpha$ is an isomorphism and we are done.
\end{proof}

\bigskip

\subsection{Relation between the Prym variety $P_\la$ and the Donagi-Prym 
$\Prym(\pi, \Sbb_\omega)$}

By Proposition \ref{evlift} there exists a morphism
$$ \widetilde{ev}_\la : \Prym(\pi, \Sbb_\omega) \lra \Jac(H; Z) $$
and because of the exact sequence \eqref{ramZY} the connected component
$\Prym(\pi, \Sbb_\omega)_0$ containing $0$ is mapped by $\widetilde{ev}_\la$
into $\Jac(Y)$.

\begin{prop} \label{imdonprym}
We have an equality of abelian subvarieties of $\Jac(Y)$
$$ \widetilde{ev}_\la \left( \Prym(\pi, \Sbb_\omega)_0 \right) 
= P_\la  \subset \Jac(Y).$$
\end{prop}

\begin{proof}
We consider the natural map of $\ZZ$-algebras
$$\phi : \bigoplus_{\omega \in \widehat{W}} \End(\Sbb_\omega)  \lra \ZZ[W]. $$
Note  that $\phi \otimes_\ZZ \QQ$ is an isomorphism of $\QQ$-algebras 
$\bigoplus_{\omega \in \widehat{W}} \End(V_\omega)  \cong \QQ[W]$. 
Taking the tensor product over $\ZZ[W]$ with the $\ZZ[W]$-module $\Jac(Z)$ induces an 
isogeny decomposition of $\Jac(Z)$ (see also \cite{D} Formula (5.3))
$$ \phi : \bigoplus_{\omega \in \widehat{W}} \Sbb_\omega \otimes_\ZZ \Prym(\pi, \Sbb_\omega)
\lra \Jac(Z) \otimes_{\ZZ[W]} \ZZ[W] = \Jac(Z).$$
Moreover the restriction of $\phi$ to $\la \cdot \ZZ  \otimes_\ZZ \Prym(\pi, \Sbb_\omega)$
is the evaluation map $ev_\la$ (section 5.1). With the notation of Lemma \ref{lem3.8},
the multiplication with the idempotent $p_\la = \frac{1}{e} S_\la \in \QQ[W]$ corresponds under 
$\phi\otimes_Z \QQ$ to the
projection onto the $\QQ$-vector 
space $\la \cdot \QQ \otimes_\QQ V^*_\omega  \subset
\bigoplus_{\omega \in \widehat{W}} \End(V_\omega)$ (see e.g. \cite{M} section 4.4).
The proposition now follows immediately from the definition of the Prym variety $P_\la$. 
\end{proof}

\bigskip

\section{The abelianization map: the \'etale case}

\subsection{Grothendieck's spectral sequence}

In this subsection we recall some facts related to Grothendieck's spectral
sequence \cite{Gr}:

\bigskip

Let $W$ be any finite 
group --- not necessarily a Weyl group, let $Z$ be a curve with a right $W$-action
with quotient $\pi : Z \ra X$
and let $A$ be an abelian algebraic group, which is also a $W$-module. Let
$\underline{A}$ denote the $W$-sheaf of abelian groups defined by $\underline{A}(U) = 
\mathrm{Mor}(U,A)$ for an open subset $U \subset Z$. Consider the following two 
left-exact functors
$$ \Gamma_Z : \{ W \text{-sheaves over} \  Z \}  \ra \{ W \text{-modules} \}, \qquad \Gamma_Z(\underline{A}) 
= \Gamma(Z,\underline{A}) = \text{global sections of }  \underline{A},$$
$$ \Gamma^W : \{ W\text{-modules} \} \ra \{ \text{abelian groups} \}, \qquad \Gamma^W(M) = M^W = 
W \text{-invariant elements of } M.$$
Consider the composite functor $\Gamma^W_Z = \Gamma^W \circ \Gamma_Z$. Its $n$-th
derived functor, which we denote by
$H^n(Z;W; \cdot)$, is computed by Grothendieck's spectral sequences 
\begin{eqnarray*}
E_2^{p,q} = H^p(W, H^q(Z,\underline{A})) \ \  & \Rightarrow  &  \ \  
E^{p+q} = H^{p+q}(Z;W;\underline{A}) \\
'E_2^{p,q} = H^p(X, R^q\pi_*^W(\underline{A})) \  \ & \Rightarrow & \ \ 
E^{p+q} = H^{p+q}(Z;W;\underline{A}) 
\end{eqnarray*}
and the associated exact sequence of low degree terms of the first spectral sequence is
\begin{equation} \label{gres0}
0 \ra E^{1,0}_2 \ra E^1 \ra E^{0,1}_2 \stackrel{c}{\ra} E^{2,0}_2 \ra E^2.
\end{equation}

In the case $W$ acts on a torus $A= T$ with $\Sbb_\omega = \Hom(T,\CC^*)$ this
exact sequence becomes using Lemma \ref{bijdp}

\begin{equation} \label{gres1}
0 \ra H^1(W,T) \ra H^1(Z;W;\underline{T}) \ra \Prym(\pi, \Sbb_\omega)
\stackrel{c}{\ra} H^2(W,T) \ra H^2(Z;W;\underline{T})
\end{equation}

One can work out the following description of the homomorphism $c$. Here we omit the details.

\begin{prop}
For any  $T$-bundle $E_T \in \Prym(\pi, \Sbb_\omega)$ the cohomology class
$c(E_T) \in H^2(W,T)$ equals the extension class \eqref{extmg} of the twisted
Mumford group $\cG^\sigma(E_T)$ of $E_T$.
\end{prop}

\begin{lem} \label{csurjetale}
If the covering $\pi: Z \ra X$ is \'etale, then $H^2(Z;W;\underline{T}) = 0$.
\end{lem}

\begin{proof}
We use the second spectral sequence $'E^{p,q}_2$. First we observe that
the sheaves $R^1\pi_*^W(\underline{T})$ and $R^2\pi_*^W(\underline{T})$
are supported on the ramification divisor of $\pi : Z \ra X$ by 
\cite{Gr} Th\'eor\`eme 5.3.1. Hence $'E^{1,1}_2 = 'E^{0,2}_2 = 0$. 

\bigskip

Next we claim that $\pi^* (\pi^W_*(\underline{T})) = \underline{T}$: let $V \subset Z$ be a sufficiently small open subset such that $g(V) \cap V = \emptyset$ for any $g \in W$
and $g \not= e$ --- here we use the analytic topology on $Z$.  Then the open 
subset $\pi^{-1} (\pi(V))$ decomposes as a disjoint union 
$\coprod_{g \in W} g(V)$. 
Now the elements of $\pi^*(\pi^W_*(\underline{T}))(V)$ correspond to $W$-equivariant 
morphisms $\Mor^W(\pi^{-1}(\pi(V)), T)$. The restriction to $V$ gives a canonical
bijection
$$ \Mor^W( \coprod_{g \in W} g(V), T) = \Mor(V,T) = \underline{T}(V),$$
which proves the claim.
\bigskip

Finally we observe that $\underline{T} \cong (\underline{\CC}^*)^n$ and since 
$H^2(Z,\underline{\CC}^*) = 0$ --- this easily follows from the
exponential exact sequence  $0 \ra \ZZ \ra \cO_Z \ra \cO^*_Z \ra 0$,
we obtain that $H^2(Z, \underline{T}) = 0$. We now conclude that
$'E^{2,0}_2 = H^2(X, \pi_*^W(\underline{T})) = 0$, since
$H^2(X, \pi_*^W(\underline{T}))$ identifies with the $W$-invariant 
subspace of $H^2(Z,\underline{T})$.
\end{proof}

\bigskip

For any $\la \in \Sbb_\omega = \Hom(T, \CC^*)$ let $H  = \Stab(\la) \subset W$ and
let $\beta_\la$ be the composite map
$$\beta_\la : H^2(W,T) \map{res_H} H^2(H,T) \map{\la} H^2(H, \CC^*).$$
If the covering $\pi: Z \ra X$ is not necessarily \'etale, 
we can say the following on the image of $c$. The
next proposition will not be used in the sequel as we will concentrate  on the
\'etale case.

\begin{prop}
We have a commutative diagram
$$
\begin{CD}
\Prym(\pi,\Sbb_\omega) @>c>> H^2(W,T) \\
@VV ev_\lambda V    @VV \beta_\lambda V \\
\Jac^H(Z) @>c'>> H^2(H, \CC^*)
\end{CD}
$$
and $\beta_\la \circ c = 0$ for any $\la \in \Sbb_\omega$.
\end{prop}

\begin{proof}
In order to prove that the diagram commutes it is sufficient to combine the 
three exact sequences \eqref{gres0} obtained in the three cases $(W,T),(H,T)$ and
$(H, \CC^*)$ --- here the first factor denotes the finite group and the 
second the abelian algebraic group $A$ on which the finite group acts. We leave the
details to the reader. The claim $\beta_\la \circ c = c' \circ ev_\la = 0$ is an
immediate consequence of Proposition \ref{evlift} and \eqref{esetale}.
\end{proof}

\bigskip

\subsection{Moduli stack of $G$-bundles}

Let $G$ be a simple and simply-connected algebraic group and let $T \subset G$
be a maximal torus, $T \subset N(T)$ the normalizer of the torus and $W = N(T)/T$ the
Weyl group of $G$. Let $\Sbb_\omega := \Hom(T, \CC^*)$ be the weight lattice 
of $G$.

\bigskip

We denote by $\cM_X(G)$ the moduli stack parametrizing principal $G$-bundles over the
curve $X$. We recall some results from \cite{LS} and \cite{S} on line bundles
over $\cM_X(G)$: let $\la \in \Sbb_\omega$ be a dominant weight and 
$\rho_\la : G \ra \mathrm{SL}(V_\la)$  the associated irreducible representation. Then
$\rho_\la$ induces a morphism of stacks
$$\widetilde{\rho}_\la : \cM_X(G) \ra \cM_X(\mathrm{SL}(m)), \qquad E_G \mapsto
E_G \times^G V_\la.$$
Here $m = \dim V_\la$.  We denote by $\DD_\la = \widetilde{\rho}^*_\la \DD$ 
the pull-back of the determinant line bundle $\DD$ over $\cM_X(\mathrm{SL}(m))$.
The next proposition is proved in \cite{LS} and \cite{S}.

\begin{prop} \label{propls}
For any simple and simply-connected algebraic group $G$ we have
\begin{enumerate}
\item There exists an ample line bundle $\cL$ over $\cM_X(G)$ such that
$$ \mathrm{Pic}(\cM_X(G)) \cong \ZZ \cdot \cL.$$
\item For any dominant weight $\la \in \Sbb_\omega$, the integer $d_\la$
defined by the relation 
$$ \DD_\la = \cL^{\otimes d_\la}$$
equals the Dynkin index of the representation $\rho_\la$ of $G$.
\end{enumerate}
\end{prop}

\bigskip

\subsection{The abelianization map $\Delta_\theta$}

>From now on we assume that the Galois covering $\pi: Z \ra X$ is \'etale with
Galois group equal to the Weyl group $W = N(T)/T$. We denote by $n \in H^2(W,T)$
the extension class of $N(T)$. Note that by Lemma \ref{csurjetale} and by \eqref{gres1} 
the homomorphism
$$ c : \Prym(\pi, \Sbb_\omega) \lra H^2(W,T)$$
is surjective. We denote by $\Prym(\pi, \Sbb_\omega)_n$ a connected component of the 
fibre of $c$ over $n$. Note that all connected components are isomorphic to
each other.

\bigskip

As in section 6.2 we consider a universal family $\cE_T$ of $T$-bundles over
$Z \times \Prym(\pi, \Sbb_\omega)_n$ --- note that $\cE_T$ is $W$-invariant --- 
and we choose an isomorphism
$$\theta: \cG^\sigma(\cE_T) \map{\sim} N(T)$$
inducing the identity on the subgroups $T$ of $\cG^\sigma(\cE_T)$ and $N(T)$ and
on the quotient $W$. The existence of $\theta$ follows from \eqref{resmumford} and
the fact that $\Prym(\pi, \Sbb_\omega)_n$ lies in the fiber of $c$ over $n$, the 
extension class of $N(T)$.

\begin{prop}
Given an isomorphism $\theta$ there exists a morphism 
$$ \Delta_\theta : \Prym(\pi, \Sbb_\omega)_n \lra \cM_X(G),$$
such that for any $E_T \in \Prym(\pi, \Sbb_\omega)_n$ the $G$-bundle 
$\Delta_\theta(E_T)$  satisfies
$$ \pi^* \Delta_\theta (E_T) = E_T \times^T G.$$
The map $\Delta_\theta$ is called the abelianization map.
\end{prop}

\begin{proof}
The existence of the morphism $\Delta_\theta$ will follow from the existence 
of a family of $N(T)$-bundles over $X$ parametrized by $\Prym(\pi, \Sbb_\omega)_n$,
or equivalently from the existence of a family of $N(T)$-bundles over $Z$ with a
$W$-linearization parametrized by $\Prym(\pi, \Sbb_\omega)_n$. But this follows
from the relative version of Proposition \ref{splitEN}, which says that
the $N(T)$-bundle $\cE_T \times^T_\theta N(T)$ admits a canonical $W$-linearization ---
here we use the isomorphism $\theta$.
\end{proof}

\begin{rem}
{\em (i): Note that $\Delta_\theta$ factorizes through $\cM_X(N(T))$. \\
(ii): A priori $\Delta_\theta$ depends on the choice of $\theta$. Two different
choices of $\theta$ differ by an element in $\Aut^0(N(T))$, the group of automorphisms
of $N(T)$ inducing the identity on $T$ and $W$. Note that $T \subset \Aut^0(N(T))$ and
that $\Aut^0(N(T))/T = H^1(W,T)$. The cohomology
groups $H^1(W,T)$ and $H^2(W,T)$ have been computed in \cite{Ma}.
}
\end{rem}

\bigskip

\subsection{Direct images of line bundles}

Consider a dominant weight $\la \in \Sbb_\omega$ and let $V_\la$ denote the
associated irreducible representation of $G$. We consider the weight space
decomposition
$$ V_\la = \bigoplus_{\mu \in \Theta(\la)} V_\la^\mu,$$
where $V_\la^\mu$ denotes the weight space of $V_\la$ associated to the
weight $\mu \in \Sbb_\omega$. The Weyl group $W$ acts on the set $\Theta(\la)$ 
of all weights $\mu$ such that $V_\la^\mu \not= \{ 0 \}$ and decomposes
it into $k$ orbit spaces 
$$ \Theta(\la) = \Theta_1 \cup \Theta_2 \cup \ldots \cup \Theta_k.$$
For each $j = 1,\ldots, k$ we choose a weight $\la_j \in \Theta_j$, i.e. $\Theta_j = 
\la_j . W$. We introduce for each $j$ the subgroup $H_j = \Stab(\la_j) \subset W$
and the quotient $\psi_j : Z/H_j = Y_j \ra X$. Let $n_j = \dim V_\la^{\la_j}$.
Note that $n_j = \dim V_\la^\mu$ for any $\mu \in \Theta_j$, that $\la_1 = \la$ and
$n_1 = 1$. 

\bigskip

Since $\pi$ is \'etale, we have $\Jac(H_j;Z) = \Jac(Y_j)$ and by Proposition \ref{evlift}
there exist morphisms
$$ \widetilde{ev}_{\la_j} : \Prym(\pi, \Sbb_\omega) \lra \Jac(Y_j). $$
For convenience of notation we introduce the product and the  map
$$\Jac( Y_\bullet ) := \Jac(Y_1) \times \Jac(Y_2) \times \cdots \times \Jac(Y_k), \qquad
\widetilde{ev}_\bullet := (\widetilde{ev}_{\la_1}, \widetilde{ev}_{\la_2}, \ldots,
\widetilde{ev}_{\la_k}),$$
and the direct image morphism
$$\psi_\bullet : \Jac( Y_\bullet) \lra \cM_X(\mathrm{GL}(m)) , \qquad 
(M_1,M_2, \ldots, M_k) \mapsto 
\bigoplus_{j=1}^k \left[ \psi_{j*}(M_j) \right]^{\oplus n_j}. $$

\begin{prop} \label{commdiagdelta}
We assume $\pi$ \'etale. Then we have a commutative diagram
$$
\begin{CD}
\Prym(\pi, \Sbb_\omega)_n @>\widetilde{ev}_\bullet>> \Jac( Y_\bullet )  \\
@VV \Delta_\theta V @VV \psi_\bullet V \\
\cM_X(G) @>\tilde{\rho}_\lambda>> \cM_X(\mathrm{GL}(m))
\end{CD}
$$
\end{prop}

\begin{proof}
Let $E_T \in \Prym(\pi, \Sbb_\omega)_n$. The rank-$m$ vector bundle 
$\tilde{\rho}_\lambda \left[ \Delta_\theta (E_T) \right] $ pulls back
under $\pi$ to the decomposable $W$-linearized vector bundle over $Z$
$$E_T \times^T V_\la =  \bigoplus_{\mu \in \Theta(\la)} E_T \times^T V_\la^\mu 
= \bigoplus_{j=1}^k \left[ \bigoplus_{\mu \in \Theta_j} E_T \times^T_\mu \CC   
\right]^{\oplus n_j}.$$
Clearly the $W$-linearization preserves the direct summands $\bigoplus_{\mu \in \Theta_j} E_T \times^T_\mu \CC $  for $j= 1, \ldots, k$ and by Lemma \ref{dirim} we have the equality
$$ \bigoplus_{\mu \in \Theta_j} E_T \times^T_\mu \CC = \pi^* ( \psi_{j*} (M_j) ),$$
as $W$-linearized vector bundles. Here $M_j = \widetilde{ev}_j(E_T)$. This proves the claim.
\end{proof}

\bigskip

\begin{rem}
{\em Note that the composite map $\psi_\bullet \circ \widetilde{ev}_\bullet = 
\tilde{\rho}_\lambda \circ \Delta_\theta$ takes values in $\cM_X(\mathrm{SL}(m))$.}
\end{rem}

\bigskip

\section{Proof of the main theorem}

For the convenience of the reader we recall the set-up of our main result.
Let $G$ be a simple and simply-connected algebraic group and $T \subset G$ a maximal torus.
Let $W := N(T)/T$ denote the Weyl group of $G$ and $\Sbb_\omega := \Hom(T, \CC^*)$ the
weight lattice. Let $\pi: Z \ra X$ be an \'etale Galois covering of smooth projective curves
with Galois group $W$. For a dominant weight $\la \in \Sbb_\omega$ we consider the Prym 
variety (see section 3.5) 
$$ P_\la \subset \Jac(Y)$$
with $Y = Z/H$ and $H = \Stab(\la)$. 
Let $\overline{K}_\la$ denote the Kanev correspondence on the curve 
$Y$. Let $L_Y$ be a line bundle over $\Jac(Y)$ representing the principal polarization. Then we can 
prove our main result.\\

\begin{thm} \label{maintheo}
Assume that the following holds:
\begin{itemize}
\item the $W$-covering $\pi : Z \ra X$ is \'etale,
\item $q_\la =  d_\la$,
\item the group $\Gamma_\la = \Sbb_\omega/ \la \cdot \ZZ[W] $ is trivial,
\item the weight $\la$ is minuscule or quasi-minuscule,
\item the homomorphism $\psi^*: \Jac(X) \ra \Jac(Y)$ is injective.
\end{itemize}
Then
\begin{enumerate}
\item there exists a line bundle $M$ over $P_\la$ such that 
$$ L_{Y}| P_\la = M^{\otimes q_\la}.$$
\item the type of the polarization $M$ over $P_\la$ equals
$$ K(M) = (\ZZ/m\ZZ)^{2g},   \qquad \text{with} \qquad m = \frac{\deg(Y/X)}
{\mathrm{gcd}(  \deg (\overline{K}_\la) -1,  \deg (Y/X) )},  $$
and $g$ denotes the genus of $X$. 
\end{enumerate}
\end{thm}

\begin{proof}
 {\bf (1)}
If $\pi$ is \'etale, Proposition \ref{evlift} gives a morphism
$$ \widetilde{ev}_\la : \Prym(\pi, \Sbb_\omega) \lra \Jac(H; Z) = \Jac(Y).$$
Moreover by Proposition \ref{imdonprym} we have 
$\widetilde{ev}_\la (\Prym(\pi, \Sbb_\omega)_0) = P_\la$. Since we have assumed that
$\Gamma_\la$ is trivial, the evaluation map $ev_\la$ is injective 
by Proposition \ref{kerevla}. Hence the map $\widetilde{ev}_\la$ is 
also injective and induces isomorphisms by restriction
$$ \widetilde{ev}_\la  :  \Prym(\pi, \Sbb_\omega)_0 \map{\sim} P_\la \subset  \Jac(Y),
\qquad   \widetilde{ev}_\la : \Prym(\pi, \Sbb_\omega)_n \map{\sim} T_\alpha (P_\la) 
\subset \Jac(Y),$$
where $T_\alpha(P_\la)$ denotes the translate by an element $\alpha \in \Jac(Y)$ of the
Prym variety $P_\la$. In order to show (1) it suffices to show that 
$L_{Y}|T_\alpha (P_\la)$ is divisible by $q_\la$.
\bigskip

We first consider the case $\la$ minuscule, i.e. $k=1$ with the notation of section 
7.4. In that case the commutative diagram of Proposition \ref{commdiagdelta} 
simplifies
$$
\begin{CD}
\Prym(\pi, \Sbb_\omega)_n @>\widetilde{ev}_\la>> T_\alpha (P_\la)   \\
@VV \Delta_\theta V @VV \psi_* V \\
\cM_X(G) @>\widetilde{\rho}_\lambda>> \cM_X(\mathrm{SL}(m)).
\end{CD}
$$
With the notation of section 7.2 we know that $\widetilde{\rho}^*_\lambda \DD = 
\cL^{\otimes d_\la}$.  Hence the line bundle 
$$  \widetilde{ev}^*_\la \left((\psi_*)^* \DD \right) = \Delta_\theta^* 
(\widetilde{\rho}^*_\lambda \DD ) = \left( \Delta_\theta^* \cL 
\right)^{\otimes d_\la}$$
is divisible by $d_\la = q_\la$. Since the first horizontal map 
$\widetilde{ev}_\la$ is an isomorphism, we deduce that $(\psi_*)^* \DD =
L_{Y}|T_\alpha (P_\la)$  is also divisible by $q_\la$.

\bigskip

In the case $\la$ is quasi-minuscule, we have $k=2$ and $\la_2 =0$. Hence  $H_2 =
\Stab(0) = W$ and $\psi_2 = \mathrm{id} : Y_2 = X \ra X$. Moreover $\widetilde{ev}_0 :
\Prym(\pi, \Sbb_\omega)_n \ra \Jac(X)$ is the constant zero map and the 
commutative diagram of Proposition \ref{commdiagdelta} simplifies
$$
\begin{CD}
\Prym(\pi, \Sbb_\omega)_n @>\widetilde{ev}_\la>> T_\alpha (P_\la)   \\
@VV \Delta_\theta V @VV \psi_\bullet V \\
\cM_X(G) @>\widetilde{\rho}_\lambda>> \cM_X(\mathrm{SL}(m)),
\end{CD}
$$
with $\psi_\bullet(M) = \psi_{1*} (M) \oplus \cO^{\oplus n_2}_X$ for $M \in 
\Jac(Y)$. We then conclude as in the case of a minuscule weight.

\bigskip

{\bf (2)} Let $S$ denote the usual Prym variety $\Prym(Y/X)$. Then we recall 
from section 3.5 that $P_\la = \im (u_\la) \subset S$ and that the endomorphism
$u_\la \in \End(S)$ satisfies the relation $u_\la^2 = q_\la u_\la$. Hence we can apply
Proposition \ref{prop2.10} (a), which says that the type of the polarization $K(M)$
is given by 
$$ K(M) = u_\la (K(L_{Y}|S)),$$
where $L_{Y}|S$ denotes the restriction of the line bundle $L_Y$ to $S$. Since 
$\psi^* : \Jac(X) \lra \Jac(Y) $ is injective, we have --- see Remark \ref{remarkprym} (ii) ---
the equalities $S = \Prym(Y/X) = \im (d-t) $ and $\psi^* \Jac(X) = \im (t)$, where $t \in \End(\Jac(Y))$
denotes endomorphism associated to the trace correspondence of $\psi: Y \ra X$. By \cite[Corollary 12.1.4]{BL} or 
Lemma \ref{lem2.5} we obtain that 
$$ K(L_{Y}|S) = K(L_{Y}| \Jac(X)) = \Jac(X) \cap S.$$
Moreover $K(L_{Y}| \Jac(X)) = \psi^* \left[ \Jac(X)_d  \right] \cong (\ZZ / d\ZZ)^{2g}$, with
$d = \deg(Y/X)$. We deduce from Corollary \ref{cor3.7}(2) that 
$u_\la t = ( \deg (\overline{K}_\la) -1 ) t$ and since $K(L_{Y}|S) \subset \im (t) = 
\psi^* \Jac(X)$ we obtain that $u_\la (K(L_{Y}|S))$ equals the image of the 
multiplication by $\deg (\overline{K}_\la) -1$ in the group  $(\ZZ / d\ZZ)^{2g}$. 
This proves the assertion (2).  
\end{proof}

\bigskip

We deduce from Lemma \ref{groupgamma} and from the two tables of section 4 the following
list of examples satisfying the 5 conditions of Theorem \ref{maintheo}. This shows the main theorem stated in the introduction.

\begin{cor}
The type $K(M) = (\ZZ/m\ZZ)^{2g}$ of the induced polarization $L_{Y| P_\la} = M^{\otimes q_\la}$
on the Prym variety $P_\la$ is given by the table.
\bigskip

\begin{center}
\begin{tabular}{|c|c||c|c|}
\hline
 Weyl group of type & weight & $q_\la = d_\la$ &  $m$ \\ \hline \hline
  $A_n$ & $\varpi_i ; \ (i,n+1) = 1$ & ${n-1 \choose i-1}$ &  $n+1$ \\ \hline 
$D_n$; \ $n$ odd & $\varpi_{n-1}, \varpi_n$  & $2^{n-3}$ & $4$ \\ \hline
$E_6$ & $\varpi_1,\varpi_6$ & $6$ & $3$ \\ \hline
$E_7$ & $\varpi_7$ & $12$ & $2$ \\ \hline
$E_8$ & $\varpi_8$ & $60$ & $1$ \\ \hline 
\end{tabular}
\end{center}
\end{cor}

\begin{proof}
We only have to check that $\psi^*: \Jac(X) \ra \Jac(Y)$ is injective in the cases mentioned in the table.
This follows from the formula $\ker \psi^* = 
\ker (W^* \ra H^*)$, where $W^*$ and $H^*$ denote the groups of characters of $W$ and $H$,
and from a straightforward case-by-case study.
\end{proof}


\bigskip

\begin{thebibliography}{999999}
\bibitem[BL]{BL} Ch. Birkenhake, H. Lange: Complex Abelian Varieties (Second Edition), Grundlehren der mathematischen Wissenschaften, Vol. 302, Springer (2004)

\bibitem[BNR]{BNR} A. Beauville, M.S. Narasimhan, S. Ramanan: Spectral covers and the generalised theta divisor,
J. Reine Angew. Math. 398 (1989), 169-179

\bibitem[Bo]{Bou} N. Bourbaki: Groupes et alg\`ebres de Lie, Chapitres 4,5 et 6, Hermann (1968)

\bibitem[Dol]{Do} I. Dolgachev: Invariant stable bundles over modular curves $X(p)$. Recent progress in algebra (Taejon/Seoul, 1997), 65-99, Contemp. Math., 224, Amer. Math. Soc.,
Providence, RI, 1999

\bibitem[Don1]{D} R. Donagi:  Decomposition of spectral covers, Journ\'ees de g\'eom\'etrie
alg\'ebrique d'Orsay, Ast\'erisque 218 (1993), 145-176

\bibitem[Don2]{D2} R. Donagi: Spectral covers, In Current Topics in complex algebraic geometry 
(Berkeley, CA, 1992/93), 65-86, Cambridge University Press, Cambridge, 1995

\bibitem[DG]{DG} R. Donagi, D. Gaitsgory: The gerbe of Higgs bundles, Transform. Groups,
No. 7 (2002), 109-153

\bibitem[Fa1]{Fa1} G. Faltings: A proof for the Verlinde formula, J. Alg. Geom. 3 (1994), 347-374

\bibitem[Fa2]{Fa2} G. Faltings: Stable $G$-bundles and projective connections, J. Alg. Geom. 2 (1993), 507-568

\bibitem[Fu]{F} W. Fulton: Intersection Theory. Erg. der Math. 2, Springer (1984).

\bibitem[Gr1]{Gr} A. Grothendieck: Sur quelques points d'alg\`ebre homologique, 
Tohoku Math. J. (2), 9 (1957), 119-221.

\bibitem[Gr2]{Gr2} A. Grothendieck: Technique de descente et th\'eor\`emes d'existence
en g\'eom\'etrie alg\'ebrique, S\'eminaire Bourbaki N. 190 (1959)

\bibitem[Ha]{hall} M. Hall: Combinatorial Theory (Second Edition), John Wiley (1986).

\bibitem[Hi1]{Hi1} N. Hitchin: The self-duality equations on a Riemann surface, Proc. London Math. Soc (3) 55
(1987), 59-126

\bibitem[Hi2]{Hi2} N. Hitchin: Stable bundles and integrable systems, Duke Math. J. 54 (1987), 91-114   

\bibitem[K1]{K} V. Kanev: Spectral curves and Prym-Tjurin varieties I, in Abelian varieties, Proc. of the Egloffstein conference 1993, de Gruyter (1995), 151-198

\bibitem[K2]{K2} V. Kanev: Principal polarizations of Prym-Tyurin varieties, Compos. Math. 64 (1987),  243-270

\bibitem[LR]{LR} H. Lange, S. Recillas: Abelian varieties with group action, J. reine angew. Math. 575 (2004), 135-155 

\bibitem[LRo]{LRo} H. Lange, A. Rojas: A Galois theoretic approach to Kanev's correspondence, Preprint (2005).

\bibitem[LS]{LS} Y. Laszlo, C. Sorger: The line bundles on the moduli of parabolic $G$-bundles over curves and their
sections, Ann. Sci. Ecole Norm. Sup. (4) 30 (1997), 
499-525

\bibitem[Ma]{Ma} M. Matthey: Normalizers of maximal tori and cohomology of Weyl groups, 
preprint available at http://igat.epfl.ch/matthey/publications/doc/weylcoho.pdf


\bibitem[Me]{M} J.-Y. M\'erindol : Vari\'et\'es de Prym d'un rev\^etement galoisien, J. reine angew. Math. 461 (1995), 49-61

\bibitem[Mu]{Mu} D. Mumford : Abelian Varieties. Tata Institute of Fundamental Research, Studies 
in Mathematics, Vol. 5, Bombay, Oxford University Press, London (1970)

\bibitem[O]{O} W. Oxbury: Spin Verlinde spaces and Prym theta functions, Proc. London Math. Soc. (3) 78 (1999), 52-76


\bibitem[Se]{Se} J.-P. Serre: Repr\'esentations lin\'eaires des groupes finis, Hermann, 
Paris (1967)

\bibitem[So1]{So1} C. Sorger: La formule de Verlinde, S\'eminaire Bourbaki, 1994/95, Exp. No. 794, Ast\'erisque 237 
(1996), 87-114



\bibitem[So2]{S} C. Sorger:  On moduli of $G$-bundles of a 
curve for exceptional $G$, Ann. Sci. Ecole Norm. Sup. (4)
32 (1999), 127-133

\bibitem[Sp]{Spr} T.A. Springer: A construction of representations of Weyl groups, Invent. Math. 44 (1978), 279-293.

 
\end{thebibliography}
\end{document}